\documentclass{elsart}
\pdfoutput=1
\usepackage{amsmath,amsfonts,amssymb,epsfig}
\newcommand{\im}{\operatorname{Im}}
\newcommand{\re}{\operatorname{Re}}

\newcommand{\spec}{\operatorname{spec}}
\newcommand{\dist}{\operatorname{dist}}
\newcommand{\R}{\mathbb{R}} %{I\!\!R}
\newcommand{\id}{I}

\newcommand{\sobol}{\mathbb{H}^1}
\newcommand{\Lint}{\mathbb{L}}
\newcommand{\C}{\mathbb{C}}
\newcommand{\Z}{\mathbb{Z}}

\newcommand{\clos}{\operatorname{cl}}
\renewcommand{\epsilon}{\varepsilon}
\renewcommand{\phi}{\varphi}

\begin{document}
\begin{frontmatter}
  \title{Longtime behavior of coupled wave equations for
    semiconductor lasers} 
  \author{Jan Sieber\thanksref{pay}}
  \thanks[pay]{This work was supported partially by the
    Sonderforschungsbereich 555 ``Komplexe Nichtlineare Prozesse'' of
    the Deutsche Forschungsgemeinschaft, and by EPSRC grant
    GR/R72020/01.}  
  \address{Dept.of Eng. Math., Queen's Building,
    University of Bristol,\\ Bristol BS8 1TR, U.K.\\ Email: Jan.Sieber@bristol.ac.uk\\
  Fax:+44 117 9546833}
\begin{abstract}
  Coupled wave equations are popular tool for investigating
  longitudinal dynamical effects in semiconductor lasers, for example,
  sensitivity to delayed optical feedback.  We study a model that
  consists of a hyperbolic linear system of partial differential
  equations with one spatial dimension, which is nonlinearly coupled
  with a slow subsystem of ordinary differential equations. We first
  prove the basic statements about the existence of solutions of the
  initial-boundary-value problem and their smooth dependence on
  initial values and parameters. Hence, the model constitutes a smooth
  infinite-dimensional dynamical system. Then we exploit the
  par\-ticu\-lar slow-fast structure of the system to construct a
  low-dimensional attracting invariant manifold for certain parameter
  constellations. The flow on this invariant manifold is described by
  a system of ordinary differential equations that is accessible to
  classical bifurcation theory and numerical tools such as AUTO.
\end{abstract}
\begin{keyword}
  laser dynamics, invariant manifold theory, strongly continuous semigroup    
\end{keyword}

% \begin{AMS}
%    78A60, 37L10, 35P10
% \end{AMS}

\end{frontmatter}

\section{Introduction}
\label{sec:intro}
Semiconductor lasers are known to be extremely sensitive to delayed
optical feedback. Even small amounts of feedback may destabilize the
laser and cause a variety of nonlinear effects. Self-pulsations,
excitability, coexistence of several stable regimes, and chaotic
behavior have been observed both in experiments and in numerical
simulations (see \cite{TA98,KL00} for general reviews). Due to their
inherent speed, semiconductor lasers are of great interest for modern
optical data transmission and telecommunication technology if these
nonlinear feedback effects can be cultivated and controlled.
Potential applications include, for example, clock recovery
\cite{PS99}, generation of pulse trains \cite{AMA95} or
high-frequency oscillations \cite{BBKSRSWH04}, and pulse reshaping
\cite{WBRH02}.

Typically, these applications utilize the laser in a non-stationary
mode, for example, to produce high-frequency oscillations or pulse
trains.  Multi-section DFB (distributed feedback) lasers allow one to
engineer these nonlinear effects by designing the longitudinal
structure of the device \cite{UBBWH04}.  If mathematical modeling is
to be helpful in guiding this difficult and expensive design process
it has to use models that are, on one hand, as accurate as possible
and, on the other hand, give insight into the nature of the observed
nonlinear phenomena. The latter is only possible by a detailed
bifurcation analysis, while only models involving partial differential
equations (PDEs) describe the effects with the necessary accuracy.

We focus in this paper on coupled wave equations with gain dispersion. This
model is a system of PDEs (one-dimensional in space), which are nonlinearly
coupled to ordinary differential equations (ODEs). It is accurate enough to
show quantitatively good correspondence with experiments and more detailed
models \cite{BBKSRSWH04,WBRH02}. We prove in this paper that the model can be
reduced to a low-dimensional system of ODEs analytically.  This makes the
model accessible to well-established and powerful numerical bifurcation
analysis tools such as \textsc{Auto} \cite{DCFKSW98}. This in turn allows us
to construct detailed and accurate numerical bifurcation diagrams for many
practically relevant situations; see \cite{WBRH02,S02} for recent results and
section~\ref{sec:out} for an illustrative example.

We achieve the central goal of our paper, the proof of the model
reduction, in three steps. First, we show that the PDE system
establishing the coupled wave model is a smooth infinite-dimensional
dynamical system, that is, it generates a semiflow that is strongly
continuous in time and smooth with respect to initial values and
parameters. Then, we exploit the particular structure of the model
which is of the form
\begin{equation}
  \label{int:geneq}
  \begin{split}
    \dot E &=\ H(n) E\\
    \dot n &= \epsilon f(n,E)
  \end{split}
\end{equation}
where the light amplitude $E\in\Lint^2([0,L];\C^4)$ is infinite-dimensional
and the effective carrier density $n\in\R^m$ is finite-dimensional. The small
parameter $\epsilon$ expresses that the carrier density $n$ operates on a much
slower time-scale than $E$.  Hence, we investigate in the second step the
growth properties of the semigroup generated by $H$ for fixed $n$, proving the
existence of spectral gaps. In the last step we construct a low-dimensional
invariant manifold for small $\epsilon$ using the general theory on the
persistence and properties of normally hyperbolic invariant manifolds for
strongly continuous semiflows in Banach spaces \cite{BLZ98,BLZ99,BLZ00}.

The paper is organized as follows. In Section~\ref{sec:twe}, we introduce the
coupled wave model as described in \cite{BWRS01} and explain the physical
meaning of all variables and parameters.  Section~\ref{sec:nt} summarizes the
results of the paper in a non-technical but precise fashion. It points out the
difficulties and the methods and theory used in the proofs. In
Section~\ref{sec:solv} we formulate the PDE system as an abstract evolution
equation in a Hilbert space and prove that it establishes a smooth
infinite-dimensional system in this setting. In this section, we consider also
inhomogeneous boundary conditions in \eqref{int:geneq} modeling optical
injection into the laser.  In Section~\ref{sec:spec} we investigate the
spectral properties of the operator $H$ for fixed $n$ and homogeneous boundary
conditions and find $n$ such that $H(n)$ has a spectral gap to the left of
the imaginary axis. Section \ref{sec:im} is concerned with the construction of
a finite-dimensional attracting invariant manifold, where we make use of the
slow-fast structure of \eqref{int:geneq} and the results of
Section~\ref{sec:solv} and Section~\ref{sec:spec}.

Finally, in Section~\ref{sec:out} we explain how the system of ODEs obtained
in Section~\ref{sec:im} can be made accessible to standard numerical
bifurcation analysis tools like AUTO. We present an illustrative and relevant
example to demonstrate the usefulness of the model reduction. Moreover, we
extend the model reduction theorem of Section~\ref{sec:im} to delay
differential equations (DDEs), which are widely used to study delayed feedback
effects in lasers \cite{TA98}.

\section{Coupled wave equations with gain dispersion}
\label{sec:twe}

The \emph{coupled wave model}, a hyperbolic system of PDEs coupled
with a system of ODEs is a well known model describing the
longitudinal effects in classical semiconductor lasers
\cite{BWRS01,TLO94}. It has been derived from Maxwell's equations for
an electro-magnetic field in a periodically modulated waveguide
\cite{BWRS01} assuming that transversal and longitudinal effects can
be separated.  In this section we introduce the corresponding system
of differential equations, explain the physical interpretation of its
coefficients and specify all physically sensible assumptions about
these coefficients.

The variable $\psi(t,z)\in\C^2$ describes the complex amplitude of the
slowly varying envelope of the optical field split into a forward and
a backward traveling wave. The variable $p(t,z)\in\C^2$ describes the
corresponding nonlinear polarization of the material. Both quantities
depend on time and the one-dimensional spatial variable $z\in[0,L]$
(the longitudinal direction within the laser). A prominent feature of
multi-section lasers is the splitting of the overall interval $[0,L]$
into sections, that is, $m$ subintervals $S_k$ that represent sections
with separate electric contacts. The other dependent variable
$n(t)\in\R^m$ describes the section-wise spatially averaged carrier
density within each section $S_k$.  In dimensionless form the
initial-boundary value problem for $\psi$, $p$, and $n$ reads as
(dropping the arguments $t$ and $z$; all coefficients are discussed below):
\begin{align}
  \partial_t
    \psi
  &=
  \begin{bmatrix}
    -1 &0\\[-1ex] \phantom{-}0&1
  \end{bmatrix}\partial_z\psi-
  \begin{bmatrix}
    0&i\\[-1ex] i&0
  \end{bmatrix}\kappa\,\psi+\bigl[(1+i\alpha_H)\,G(n)-d\bigr]\,\psi
  +\rho(n)\,\bigl[p-\psi\bigr]
\label{twe:twe}\\[0.8ex]
  \partial_t  p&=[i\Omega_r(n)-\Gamma(n)]\,p+
  \Gamma(n)\,\psi \label{twe:pol}%\\[0.8ex]
  \intertext{and for $k=1\ldots m$}
  \dot n_k&=
  I_k-\frac{n_k}{\tau_k}-\frac{P}{l_k}
   \int\limits_{z_k}^{z_{k+1}}\left[
     \begin{array}[c]{cc}
\psi\\ p
\end{array}
\right]^*
   \left[\begin{array}[c]{ccc}
     G_k(n_k)-\rho_k(n_k) &\ & \frac{1}{2}\rho_k(n_k) \\
    \frac{1}{2}\rho_k(n_k) &\ & 0    
   \end{array}\right]
 \left[
   \begin{array}[c]{c}
     \psi\\ p
   \end{array}\right]\mbox{d}z\label{twe:carrier}\\[1ex]
 &=:f_k(n_k,(\psi,p))\mbox{.}\nonumber
\end{align}
 System
\eqref{twe:twe}--\eqref{twe:carrier} is subject to the inhomogeneous
boundary conditions for $\psi$
\begin{equation}
  \label{twe:bcinh}
  \psi_1(t,0)=r_0\psi_2(t,0)+\alpha(t)\mbox{,\ }
  \psi_2(t,L)=r_L\psi_1(t,L)
\end{equation}
and the initial conditions 
\begin{equation}
  \label{twe:inicond}
  \psi(0,z)=\psi^0(z)\mbox{,\ } p(0,z)=p^0(z)\mbox{,\ } n(0)=n^0\mbox{.}
\end{equation}
The Hermitian transpose of a $\C^4$ vector  is denoted by
$[\cdot]^*$ in \eqref{twe:carrier}. 
The length of the laser is $L$. We denote the length of
subinterval $S_k$ by $l_k$ and its starting point by $z_k$ (for
$k=1\ldots m$). We scale $z$ such that $l_1=z_2=1$ and denote
$z_{m+1}=L$. Thus, $S_k=[z_k,z_{k+1}]$. All coefficients of
\eqref{twe:twe}, \eqref{twe:pol} are spatially constant in each
subinterval $S_k$ and depend only on $n_k$, that
is, if $z\in S_k$,
\begin{xxalignat}{3}
  \kappa(z)&=\ \kappa_k\in\R\mbox{,\quad} &
  d(z)&=d_k\in\C\mbox{,\quad} &
  \alpha_H(z)&= \alpha_{H,k}\in\R\mbox{,}\\
  G(n,z)&=G_k(n_k)\in\R\mbox{,\quad} & \rho(n,z)&=
  \rho_k(n_k)\in\R\mbox{,\quad} &
  \Omega_r(n,z)&= \Omega_{r,k}(n_k)\in\R\mbox{,}\\
  \Gamma(n,z)&=\ \Gamma_k(n_k)\in\R\mbox{.}
\end{xxalignat}
The permissible range of $n_k$ is the interval
$(\underline{n},\infty)$ where, typically, $\underline{n}=-\infty$ or
$\underline{n}=0$. We assume that the functions $\rho$, $\Omega$, and
$\Gamma$ are Lipschitz continuous, and smooth in
$(\underline{n},\infty)$. Moreover, $\Omega$ and $\Gamma$ are bounded
globally. The function
$G_k:(\underline{n},\infty)\to\R$ is a smooth strictly monotone
increasing function satisfying $G_k(1)=0$, $G_k'(1)>0$,
$\lim_{\nu\searrow\underline{n}} G_k(\nu)=-\infty$, and
$\lim_{\nu\to\infty} G_k(\nu)=\infty$. Typical models for $G_k$ are
\begin{equation}\label{gmodel}
\begin{aligned}
  G_k(\nu)&=\ \tilde g_k\log\nu &
    \mbox{(}\underline{n}&=\ 0\mbox{) or}\\ 
  G_k(\nu)&=\ \tilde g_k\cdot(\nu-1) &
  \mbox{(}\underline{n}&=\ -\infty\mbox{)}
\end{aligned}
\end{equation}
where $\tilde g_k=G_k'(1)>0$. Physically sensible assumptions on the
parameters throughout the paper are $\re d_k>0$, $\Gamma_k(n_k)>1$,
$0<I_k\ll1$, $\tau_k\gg1$.  Furthermore, $\rho_k(n_k)$ is bounded for
$n_k\searrow\underline{n}$, $\rho_k(n_k)\geq0$ for $n_k\geq1$ and
$\rho_k(n_k)/G_k(n_k)\to0$ for $n_k\to\infty$. In \eqref{twe:bcinh}
$r_0$ and $r_L\in\C$ satisfy $|r_0|<1$ and $|r_L|<1$. All other
coefficients (as well as their Lipschitz constants if they depend on
$n$) are of order $1$. The constant $P$ determines the
scaling of $(\psi,p)$. It can be chosen arbitrarily. We denote the
right-hand-side of the carrier density equation \eqref{twe:carrier}
for $n_k$ by $f_k(n_k,(\psi,p))$ observing that $f_k$ is Hermitian in
$(\psi,p)$. Collecting $\psi$ and $p$ into one $\C^4$ vector
$E=(\psi,p)^T$ system~\eqref{twe:twe}--\eqref{twe:carrier} assumes the
form of \eqref{int:geneq} discussed in the introduction. In
particular, Equations~\eqref{twe:twe}, \eqref{twe:pol} are linear in
$\psi$ and $p$ if the inhomogeneity $\alpha(t)$ equals zero.

\paragraph*{Physical background and motivation of the model}
Equation \eqref{twe:twe} describes the transport of the two complex
profiles $\psi_1$ (forward) and $\psi_2$ (backward) in a waveguide
from $z=0$ to $z=L$. The boundary conditions \eqref{twe:bcinh}
describe the reflection at the facettes at $z=0$ and $z=L$ of the
laser with complex reflectivities of modulus less than $1$ and,
possibly, a complex optical input $\alpha(t)$. The case of
discontiuous optical input $\alpha(t)$ is of interest in applications
(for example, square-wave signals, or noisy signals).

The number of sections, $m$, is typically small. For example, the prototype
device studied in \cite{WBRH02,S02} has $m=2$. The coupling between the wave
profiles $\psi_1$ and $\psi_2$, expressed by a non-zero $\kappa$ in
\eqref{twe:twe}, is achieved by a small-scale spatial modulation of the
refractive index of the waveguide. Equation~\eqref{twe:twe} considers only the
spatially averaged coupling effect, called distributed feedback, which affords
the name DFB laser.

The function $G_k$, called gain, expresses the dependence of the
amplification of $\psi$ (stimulated emission of photons) on the
carrier density $n_k$ in each section. The gain $G_k(n_k)$ is positive
if $n_k>1$ and negative if $n_k<1$.  The real part of $d_k$ is
positive expressing the wave\-guide losses. A nonzero coefficient
$\alpha_{H,k}$ (called the Henry factor, typically
$\alpha_{k,H}\in(0,10)$) expresses the fact that an increase of the
carrier density $n_k$ not only increases the amplification in the
waveguide but also changes its refractive index. This effect is
regarded as one of the main reasons behind the extreme sensitivity of
semiconductor lasers with respect to optical feedback \cite{KL00}.

Taking into account the polarization $p$ introduces the effect of
\emph{gain dispersion}, which means that the light amplification
depends on the frequency of the light, prefering frequencies close to
$\Omega_r$. This effect is rather weak in long semiconductor lasers
(which means $\Gamma\gg1$) but it breaks the symmetry
$[\lambda+\re d-G(n),E(z)]\to[G(n)-\re d-\bar\lambda,\bar
E(L-z)]$ of the eigenvalue problem $\lambda E=H(n)E$ corresponding to
the fast subsystem in \eqref{int:geneq} in the case of a single
section ($m=1$). This symmetry breaking has a significant effect on
the observed dynamics also in the multi-section case \cite{BWRS01}.
In combination with the non-zero coupling $\kappa$ the presence of
a positive $\rho_k(n_k)$ for at least
one section $S_k$ and a finite $\Gamma_k(n_k)$ guarantee the existence
of a spectral gap for the linear operator $H(n)$ in \eqref{int:geneq},
which will be needed for the construction of the invariant manifold.

The dynamics of the carrier density (\eqref{twe:carrier} in section
$S_k$) is governed by three terms: the constant current input $I_k$, a
decay with rate $\tau_k^{-1}$, and the \emph{stimulated recombination},
spatially averaged over $S_k$. The stimulated recombination is
Hermitian in $E=(\psi,p)^T$. The coefficients $G_k(n_k)$ and
$\rho_k(n_k)$ are such that the $n_k$-dependent matrix in the
integrand is positive definite for $n_k\to\infty$ and negative
definite for $n_k\searrow\underline{n}$. The smallness of the
parameters $I_k$ and $\tau_k^{-1}$ in the appropriate
non-dimensionalization and the choice of a small $P$ makes $\dot n_k$
small. That is, $n$ is a slow variable compared to $E=(\psi,p)^T$
\cite{TA98}.

\section{Non-technical overview}
\label{sec:nt}
In this section we state the main results of the paper in a
non-technical but precise manner and summarize the methods used in the
proofs of these results. We have split this section into four parts.
First we show that system \eqref{twe:twe}--\eqref{twe:carrier}
generates a smooth infinite-dimensional dynamical system. Then we
introduce a small parameter. In the next step we investigate the dynamics of the
(linear) infinite-dimensional fast subsystem, and finally we construct
a low-dimensional attracting invariant manifold.

\subsection{Existence theory}\label{sec:nt:solv}
In a first step we investigate in which sense system
\eqref{twe:twe}--\eqref{twe:carrier} generates a semiflow depending
smoothly on its initial values and all parameters; for details see
section~\ref{sec:solv}. Our aim is to write
\eqref{twe:twe}--\eqref{twe:carrier} as an abstract evolution equation
in the form
\begin{displaymath}
  \frac{d}{dt}u=Au+g(u)
\end{displaymath}
in a Hilbert space $V$ where $A$ is a linear differential operator
that generates a strongly continuous semigroup $S(t)$ and $g$ is
smooth in $V$. A natural space for the variables $\psi$ and $p$ is
$\Lint^2([0,L];\C^2)$, such that $V$ could be
$\Lint^2([0,L];\C^2)\times \Lint^2([0,L];\C^2)\times\R^m$ for the
variable $u=(\psi,p,n)$. However, the inhomogeneity $\alpha$ in the
boundary condition \eqref{twe:bcinh} does not fit into this framework.
Common approaches such as boundary homogenization (used in
\cite{PS99}) % or appending $\alpha$ as a simple auxiliary variable
% (done in \cite{F99})
require a high degree of regularity of $\alpha$
in time, which is unnatural as the laser still works with discontinuous
input.  To avoid the introduction of the concept of weakly mild solutions (as was
done in \cite{JR02}), we introduce the auxiliary space-dependent
variable $a(t,x)$ ($x\in[0,\infty)$) satisfying the equation
\begin{equation}\label{nt:inh}
  \partial_t a(t,x)=\partial_x a(t,x)
\end{equation}
and change the boundary condition for $z=0$ in \eqref{twe:bcinh} into 
\begin{math}
  \psi_1(t,0)=r_0\psi_2(t,0)+a(t,0)\mbox{.}
\end{math}
One may think of an infinitely long fibre $[0,\infty)$ storing all
future optical inputs and transporting them to the laser facet $z=x=0$
by the transport equation \eqref{nt:inh}. If we choose
$a(0,x)=\alpha(x)$ as initial value for $a$ then the value of $a$ at the
boundary $x=0$ at time $t$ is $\alpha(t)$. In this way, the formerly
inhomogeneous boundary condition becomes linear in the variables
$\psi$ and $a$ requiring no regularity for $a$. To keep the space $V$
a Hilbert space, we choose a weighted $\Lint^2$ norm for $a$ that
contains $\Lint^\infty$, that is,
$\|a(t,\cdot)\|^2=\int_0^\infty|a(t,x)|^2(1+x^2)^\eta\,dx$ with
$\eta<-1/2$.

With this modification we can work within the framework of the theory
of strongly continuous semigroups \cite{P83}. The variable $u$ has the
components $(\psi,p,n,a)\in V=\Lint^2([0,L];\C^2)\times
\Lint^2([0,L];\C^2)\times\R^m\times \Lint^2_\eta([0,\infty);\C)$. We
have a certain freedom how to choose the splitting of the
right-hand-side between $A$ and $g$. We keep $A$ as simple as possible,
including only the unbounded terms
\begin{displaymath}
  A\,\left[(\psi,p,n,a)^{\,T}\right]
%%   \begin{bmatrix}
%%     \psi\\ p\\ n\\ a
%%   \end{bmatrix}
:= \left([-\partial_z\psi,\partial_z\psi],0,0,\partial_xa\right)^{\,T}
%% \begin{bmatrix}
%%   \begin{bmatrix}
%%     -\partial_z\psi_1\\
%%     \phantom{-}\partial_z\psi_2
%%   \end{bmatrix}\\
%%   0\\ 0\\ \partial_x a
%% \end{bmatrix}
\mbox{.}
\end{displaymath}
In this way, it is easy to prove that $A$ generates a strongly
continuous semigroup $S(t)$ by constructing $S$ explicitly. The
nonlinearity $g$ is smooth because it is a superposition operator of
smooth coefficient functions, and all components either depend only
linearly on the infinite-dimensional components $\psi$ and $p$, or map
into $\R^m$. Then, the existence of a semiflow $S(t;u)$ that is
strongly continuous in $t$ and smooth with respect to $u$ and
parameters follows from an a-priori estimate. This a-priori estimate
is more subtle than the one in \cite{PS99}. It uses the fact that
there is dissipation (for $\re d_k>0$) and that the same functions
$G_k$ and $\rho_k$ appear on the right-hand-side of \eqref{twe:twe}
and of \eqref{twe:carrier} but with opposing signs.  Due to this fact
we can show that the function
\begin{displaymath}
  \frac{P}{2}\|\psi(t)\|^2+\sum_{k=1}^m l_k (n_k(t)-n_*)
\end{displaymath}
remains non-negative for sufficiently small $n_*$ and, hence, bounded,
giving rise to a bounded invariant ball in $V$.

\subsection{Introduction of a small parameter}\label{sec:nt:epsilon}
For all results about the long-time behavior of system
\eqref{twe:twe}--\eqref{twe:carrier} we restrict ourselves to
autonomous boundary conditions for $\psi$, that is,
\begin{equation}
  \label{nt:bcond}
   \psi_1(t,0)=r_0\psi_2(t,0)\mbox{,\quad }
  \psi_2(t,L)=r_L\psi_1(t,L)\mbox{.}
\end{equation}
The inhomogeneous case is an open question for future work. However,
un\-der\-stan\-ding the dynamics of the autonomous laser is not only an
intermediate step but an important goal in itself since many
experiments and simulations focus on this case; see, for example,
\cite{BBKSRSWH04} for further references.

We exploit that  $I_k$ and $\tau_k^{-1}$ in \eqref{twe:carrier} are
approximately two orders of magnitude smaller than $1$ by introducing   a small
parameter $\epsilon$. We set $P=\epsilon$  in \eqref{twe:carrier}  such
that the set of carrier density equations \eqref{twe:carrier} reads as
\begin{equation}
  \label{spec:carrier}
  \frac{d}{dt}n=f\left(n,E\right)=\epsilon F(n,E)
\end{equation}
where all coefficients of $F=f/\epsilon$ are of order $1$ in all $m$
components of $F$, and $E=(\psi,p)^T$.  The space dependent subsystem
\eqref{twe:twe}, \eqref{twe:pol} is linear in $E$
\begin{equation}
  \label{spec:linsub}
  \frac{d}{dt} E=H(n)E
\end{equation}
where $H(n)$ has the form
\begin{equation}
  \label{spec:h}
H(n)=  \begin{bmatrix}
    \begin{bmatrix}
      -\partial_z +\beta(n) &-i\kappa \\
      -i\kappa & \partial_z +\beta(n)
    \end{bmatrix} & \rho(n) \\[1em]
    \Gamma(n) & i\Omega_r(n)-\Gamma(n)
  \end{bmatrix}
\end{equation}
and $\beta(n)=(1+i\alpha_H)G(n)-d-\rho(n)$. For fixed $n$, $H(n)$
acts from
\begin{multline*}
  Y:=\{ (\psi,p)\in \sobol([0,L];\C^2)\times \Lint^2([0,L];\C^2):\\
 \psi_1(0)=r_0\psi_2(0)\mbox{,\ }\psi_2(L)=r_L\psi_1(L)\} 
\end{multline*}
into $X=\Lint^2([0,L];\C^2)\times\Lint^2([0,L];\C^2)$. The operator
$H(n)$ generates a $C_0$ semigroup $T(n;t):X\to X$. The coefficients
$\kappa$, and, for each $n\in(\underline{n},\infty)^m$, $\beta(n)$,
$\Omega_r(n)$, $\Gamma(n)$ and $\rho(n)$ are linear operators in
$\Lint^2([0,L];\C^2)$ defined by the corresponding coefficients in
\eqref{twe:twe}, \eqref{twe:pol}. The maps $\beta$, $\rho$, $\Gamma$,
and $\Omega_r:\R^m\to\mathcal{L}(\Lint^2([0,L];\C^2))$ are smooth.

Although $\epsilon$ is not directly accessible, we treat it as a
parameter and consider the limit $\epsilon\to 0$ while keeping $F$
fixed. At $\epsilon=0$, the carrier density $n$ is constant. It enters
the linear subsystem \eqref{spec:linsub} as a parameter. Thus, the
spectral properties of $H(n)$ determine the long-time behavior of
system \eqref{spec:carrier}, \eqref{spec:linsub} for $\epsilon=0$.
Section~\ref{sec:nt:im} (and Section~\ref{sec:im} in detail) will
study in which sense this fact can be exploited to find exponentially
attracting invariant manifolds also for small $\epsilon$.

\subsection{Spectral properties of $H(n)$}\label{sec:nt:spec}
The goal of Section~\ref{sec:spec} is to show that there exist carrier
densities $n$ such that $H(n)$ (and, correspondingly, $T(n;t)$) has a
spectral gap. More precisely, we investigate if there exists a
$n\in(\underline{n},\infty)^m$,  and a splitting of
$X=X_c(n)\oplus X_s(n)$ into two $H(n)$-invariant subspaces such that
$X_c(n)$ is finite-dimensional, the spectrum of $H(n)\vert_{X_c(n)}$
is on the imaginary axis, and the semigroup $T(n;\cdot)\vert_{X_s(n)}$
satisfies 
\begin{equation}\label{nt:specgap}
  \|T(n;t)\vert_{X_s}\|\leq M e^{-\xi t}
\end{equation}
for some $\xi>0$ and $M\geq1$ and all $t\geq0$. Carrier densities $n$ for
which $H(n)$ has this splitting are called \emph{critical}. Generically they
form submanifolds (with boundary) of $\R^m$ of dimension $m-q$ where
$q=\dim_\C X_c(n)$.

The strongly continuous semigroup $T(n;t)$ is a compact perturbation of the
semigroup $T_0(n;t)$ generated by the diagonal part $H_0(n)$ of $H(n)$
\cite{NRL86}. The essential spectral radius $\exp(tR_\infty(n))$ of $T_0$ can
be found analytically.  To prove the existence of a critical carrier density
$n$ with spectral gap we have to find a $n\in(\underline{n},\infty)^m$ such
that $R_\infty(n)<0$ but the characteristic function $h(n;\cdot)$ of $H(n)$
(which is known analytically) has roots on the imaginary axis. This is
possible if one of the following two conditions is satisfied
\begin{enumerate}
\item The coupling $\kappa_k$ is nonzero for at least
  one $k\in\{1,\ldots,m\}$, and $m\geq2$.
\item The reflectivities satisfy $r_0r_L\neq0$, and $\rho_k(\nu)>0$
  for all $\nu>1$ for at least one $k\in\{1,\ldots,m\}$.
\end{enumerate}
These two conditions state precisely in which sense the coupling
$\kappa_k$ and the gain dispersion $\rho_k(n)$ guarantee the spectral gap
of $H(n)$ (as mentioned already in Section~\ref{sec:twe}).

\subsection{Existence of a low-dimensional invariant manifold}\label{sec:nt:im}

Let $\mathcal{K}\subset \R^m$ be a manifold of critical carrier densities with
uniform spectral gap. The dimension of the critical subspace $X_c(n)$ is
constant, say, $q\geq1$. For $\epsilon=0$ and a sufficiently small neighborhood
$U$ of $\mathcal{K}$ the set
\begin{displaymath}
  \mathcal{C}_0:=\{(E,n): n\in U, E\in X_c(n), \|E\|\leq R\}
\end{displaymath}
is invariant and exponentially attracting in its $E$-component.
Theorem~\ref{im:thm} in Section~\ref{sec:im}, the main nonlinear result,
proves that this invariant manifold persists also for small $\epsilon>0$. For
$k\geq2$, and sufficiently small $\epsilon>0$ and $U\supset\mathcal{K}$,
system~\eqref{spec:carrier},\eqref{spec:linsub} has a $C^k$ manifold
$\mathcal{C}$ that is invariant and exponentially attracting relative to
$\{(E,n): \|E\|\leq R, n\in U\}$.  It can be described as a $C^1$-small graph
over $\mathcal{C}_0$.  The semiflow restricted to $\mathcal{C}$ is governed by
the ODE in $\C^q\times\R^m$
\begin{equation}
  \label{nt:flowc}
  \begin{split}
    \frac{d}{dt} E_c &=\ \left[H_c(n)+\epsilon
      a_1(E_c,n,\epsilon) +\epsilon^2a_2(E_c,n,\epsilon)
      \nu(E_c,n,\epsilon)\right] E_c\\ 
    \frac{d}{dt} n_{\phantom{k}} &=\ \epsilon F(n,[B(n)+
    \epsilon\nu(E_c,n,\epsilon)]E_c)
  \end{split}
\end{equation}
where $H_c(n)=H(n)\vert_{X_c(n)}$, $B(n)$ is a basis of $X_c(n)$, $\nu\in
C^{k-2}$, and
\begin{displaymath}
  \begin{split}
    a_1(E_c,n,\epsilon) &=\
    -B(n)^{-1}P_c(n)\partial_nB(n)F(n,[B(n)+
    \epsilon\nu(E_c,n,\epsilon)]E_c)\\
    a_2(E_c,n,\epsilon) &=\
    B(n)^{-1}\partial_nP_c(n)F(n,[B(n)+
    \epsilon\nu(E_c,n,\epsilon)]E_c)[\id-P_c(n)]
  \end{split}
\end{displaymath}
where $P_c(n)$ is the spectral projection onto $X_c(n)$ for $H(n)$.

Theorem~\ref{im:thm} is an application of the general theory about persistence
of normally hyperbolic invariant manifolds of semiflows under $C^1$ small
perturbations \cite{BLZ98,BLZ99,BLZ00}. The unperturbed
invariant manifold $\mathcal{C}_0$ is finite-dimensional and exponentially
stable. The proof describes in detail the appropriate cut-off modification of
the system outside of the region of interest to make the unperturbed invariant
manifold compact.  Then it connects the results of the previous sections to
guarantee the $C^1$-smallness of the perturbation and the normal hyperbolicity
of $\mathcal{C}_0$.
A similar model reduction result has been proven in \cite{T01} for ODEs of the
structure \eqref{int:geneq} using Fenichel's Theorem for singularly perturbed
systems of ODEs \cite{F79}. Since Fenichel's Theorem is not available for
infinite-dimensional systems, we have to adapt the proof in \cite{F79} to our
case starting from the general results in \cite{BLZ98,BLZ99,BLZ00} about
invariant manifolds of semiflows in Banach spaces. In particular, we apply the
cut-off modifications done in \cite {F79} only to the finite-dimensional
components $E_c$ and $n$.  Moreover, we adapt the modifications such that the
invariant manifold for $\epsilon=0$ is compact without boundary as required by
the theorems in \cite{BLZ98}.

%Apparently, the graph of the manifold $\mathcal{ C}$ itself enters
% of the form
% \begin{equation}
%   \label{nt:flowapp}
%   \begin{split}
%     \frac{d}{dt} E_c &=\ \left[H_c(n)+
%       a_1(E_c,n,\epsilon)\right] E_c\\ 
%     \frac{d}{dt} n_k &=\ \epsilon \left(F_k(n_k)-g_k(n_k)[B(n)E_c,B(n)E_c]\right)
%   \end{split}
% \end{equation}
% : 
% \begin{displaymath}
%   \begin{split}
%     H_c(n)&=B^{-1}(n)P_c(n) H(n) B(n)\\
%     a_1(E_c,n,\epsilon)& =-B^{-1}(n)P_c(n)\partial_n B(n)\, \dot n.    
%   \end{split}
% \end{displaymath}
%% Truncating all terms of order
%% $O(\epsilon^2)$ in \eqref{nt:flowc} gives rise to a system of ODEs in
%% $\C^q\times\R^m$ where all terms in the right-hand-side can be
%% expressed analytically as functions of the eigenvalues of $H$.
%% The truncated system \eqref{nt:flowc} is called the \emph{mode
%%   approximation}. It is an implicit system of ODEs because the
%% eigenvalues of $H$ are given only implicitly as roots of the
%% characteristic function $h$ of $H$. The dimension of
%% \eqref{nt:flowc} is typically low: $q$ is often either $1$ or $2$.
%% The consideration of mode approximations has proven to be extremely
%% useful for numerical and analytical investigations of longitudinal
%% effects in multi-section semiconductor lasers; see for example
%% \cite{BRS98}, \cite{S02}, \cite{WBRH02}.

\section{Existence of a smooth semiflow}
\label{sec:solv}
In this section, we treat the inhomogeneous initial-boundary value
problem \eqref{twe:twe}-\eqref{twe:bcinh} as an autonomous nonlinear
evolution equation
\begin{equation}
  \label{twe:evolsys}
  \frac{d}{dt}u(t)=A u(t)+g(u(t))\mbox{,\hspace*{1em}} u(0)=u_0
\end{equation}
where $u(t)$ is an element of a Hilbert space $V$, $A$ is a generator
of a $C_0$ semigroup $S(t)$, and $g:U\subseteq V\to V$ is smooth and
locally Lipschitz continuous in an open set $U\subseteq V$.  The
inhomogeneity in \eqref{twe:bcinh} is included in \eqref{twe:evolsys}
as a component of $u$.
%\subsection{Notation}
%\label{sec:spaces}

The Hilbert space $V$ is defined by
\begin{equation}
  \label{twe:Xdef}
  V:=\Lint^2([0,L];\C^2)\times\Lint^2([0,L];\C^2)\times 
  \R^{m}\times\Lint^2_\eta([0,\infty);\C)
\end{equation}
where $\Lint^2_\eta([0,\infty);\C)$ is the space of weighted
square integrable functions. The scalar product of
$\Lint^2_\eta([0,\infty);\C)$ is defined by
\begin{displaymath}
  (v,w)_\eta:=\re \int_0^\infty \bar v(x)\cdot w(x) (1+x^2)^\eta dx\mbox{.}
\end{displaymath}
We choose $\eta<-1/2$ such that the space $\Lint^\infty([0,\infty);\C)$ is
continuously embedded in $\Lint^2_\eta([0,\infty);\C)$.  The complex
plane is treated as two-dimensional real plane in the definition of
the vector space $V$ such that the standard $\Lint^2$ scalar product
$(\cdot,\cdot)_V$ of $V$ is differentiable.  The corresponding
components of a vector $v\in V$ are denoted by
\begin{displaymath}
v=(\psi, p, n, a).
\end{displaymath}
Here, $\psi$ and $p$ have two complex
components and $n\in\R^m$.  The spatial variable in $\psi$ and $p$ is
denoted by $z\in[0,L]$, whereas the spatial variable in $a$ is denoted
by $x\in[0,\infty)$.  The Hilbert space $\sobol_\eta([0,\infty);\C)$
equipped with the scalar product
\begin{displaymath}
  (v,w)_{1,\eta}:=(v,w)_\eta+(\partial_x v,\partial_x w)_\eta
\end{displaymath}
is densely and continuously embedded in
$\Lint^2_\eta([0,\infty);\C)$. Moreover, its elements are 
continuous \cite{T78}. Consequently, the Hilbert spaces
\begin{eqnarray*}
  W&:=& \sobol([0,L];\C^2) \times\Lint^2([0,L];\C^2)
  \times\R^m \times\sobol_\eta([0,\infty);\C)\mbox{,\ and}\\ 
  W_{\mathrm{BC}}&:=&\{(\psi,p,n,a)\in W:
  \psi_1(0)=r_0\psi_2(0)+a(0), \psi_2(L)=r_L\psi_1(L)\}
\end{eqnarray*}
are densely and continuously embedded in $V$.  The linear functionals
$\psi_1(0)-r_0\psi_2(0)-a(0)$ and $\psi_2(L)-r_L\psi_1(L)$ are
continuous from $W\to\R$. We define the linear operator
$A:W_{\mathrm{BC}}\to V$ by
\begin{equation}\label{twe:adef}
  A 
  \begin{bmatrix}
    \psi\\[-1ex] p\\[-1ex] n\\[-1ex] a
  \end{bmatrix}:= 
\begin{bmatrix}
  \begin{bmatrix}
    -\partial_z &0\\[-1ex]
    \phantom{-}0 & \partial_z
  \end{bmatrix}\,\psi\\[-1ex]
  0\\[-1ex] 0\\[-1ex] \partial_x a
\end{bmatrix}\mbox{.}
\end{equation}
The definition of $A$ and $W_\mathrm{BC}$ treat the inhomogeneity
$\alpha$ in the boundary condition \eqref{twe:bcinh} as the boundary
value of the variable $a$ at $0$.  We define the open set $U\subseteq
V$ by
\begin{displaymath}
  U:=\{(\psi,p,n,a)\in V: n_k>\underline{n}\mbox{\ for\ }k=1\ldots
  m\}\mbox{,}
\end{displaymath}
and the nonlinear function $g:U\to V$ by
\begin{equation}\label{twe:gdef}
  g(\psi,p,n,a)=
  \begin{pmatrix}
    \left[(1+i\alpha_H)G(n)-d\right]\psi -i\kappa
    \begin{bmatrix}
      0&1 \\[-1.5ex] 1& 0
    \end{bmatrix}
\psi + \rho(n)[p-\psi]\\[1em]
    (i\Omega_r(n)-\Gamma(n)) p+ \Gamma(n)\psi\\[0.5em]
      \left(f_k(n_k,(\psi,p))\right)_{k=1}^m\\
    0
  \end{pmatrix}\mbox{.}
\end{equation}
The corresponding coefficients of \eqref{twe:twe}--\eqref{twe:carrier}
define the smooth maps
 $G,\rho,\Omega_r,\Gamma:(\underline{n},\infty)^m\to \mathcal{L}(\Lint^2([0,L];\C^2))$.
The function $g$ is continuously differentiable to any order with
respect to all arguments and its Frechet derivative is bounded in any
closed bounded ball $B\subset U$. Thus, \eqref{twe:evolsys}
with the right-hand-side defined by \eqref{twe:adef}, \eqref{twe:gdef}
fits into the framework of the theory of $C_0$ semigroups \cite{P83}.

The inhomogeneous initial-boundary value problem
\eqref{twe:twe}-\eqref{twe:inicond} and the au\-ton\-omous evolution
system \eqref{twe:evolsys} are equivalent in the following sense:
Suppose $\alpha\in\sobol([0,T);\C)$ in \eqref{twe:bcinh}.  Let
$u=(\psi,p,n,a)$ be a classical solution of \eqref{twe:evolsys},
satisfying \eqref{twe:evolsys} in the $\Lint^2$ sense.
Then $u$ satisfies \eqref{twe:twe}-\eqref{twe:pol}, and
\eqref{twe:inicond} in $\Lint^2$ and \eqref{twe:carrier},
\eqref{twe:bcinh} for each $t\in[0,T]$ if and only if
$a^0\vert_{[0,T]}=\alpha$.  On the other hand, assume that
$(\psi,p,n)$ satisfies \eqref{twe:twe}-\eqref{twe:pol}, and
\eqref{twe:inicond} in $\Lint^2$ and \eqref{twe:carrier},
\eqref{twe:bcinh} for each $t\in[0,T]$. Then, choosing
$a^0\in\sobol_\eta([0,\infty);\C)$ such that $a^0\vert_{[0,T]}=\alpha$,
we obtain that $u(t)=(\psi(t),p(t),n(t),a^0(t+\cdot))$ is a classical
solution of \eqref{twe:evolsys} in $[0,T]$.

We prove in Lemma \ref{twe:generator} that $A$ generates a $C_0$
semigroup $S(t)$ in $V$.  Mild solutions of \eqref{twe:evolsys},
satisfying the variation of constants formula in $V$
\begin{displaymath}
  u(t)=S(t)u_0+\int_0^t S(t-s)g(u(s)) ds\mbox{,}
\end{displaymath}
are a sensible generalization of the classical solution concept of
\eqref{twe:twe}-\eqref{twe:bcinh} to boundary conditions including
discontinuous inputs $\alpha\in\Lint^2_\eta([0,\infty);\C)$.
%\subsection{Existence of a semiflow}
%\label{sec:semi}
%In order to prove uniqueness and global existence of solutions of
%\eqref{twe:evolsys}, we apply the theory of strongly continuous
%semigroups  \cite{P83}.
\begin{lem}\label{twe:generator}
  The operator $A:W_\mathrm{BC}\subset V\to V$ generates a $C_0$
  semigroup $S(t)$ of bounded operators in $V$.
\end{lem}
\begin{pf}%\ \\
  We specify the $C_0$ semigroup $S(t)$ explicitly. Denote the components of
  $S(t)((\psi_1^0, \psi_2^0),p^0,n^0,a^0)$ by $((\psi_1(t,z),
  \psi_2(t,z)), p(t,z), n(t), a(t,x))$ for $z\in[0,L]$,
  $x\in[0,\infty)$, and let $t\leq L$.
  \begin{displaymath}
  \begin{split}
    \psi_1(t,z)&=\ 
    \begin{cases}
      \psi_1^0(z-t) & \mbox{for $z>t$}\\
      r_0\psi_2^0(t-z)+a^0(t-z) & \mbox{for $z\leq t$}
    \end{cases}\\
    \psi_2(t,z)&=\ 
    \begin{cases}
      \psi_2^0(z+t) & \mbox{for $z<L-t$}\\
      r_L\psi_1^0(2L-t-z) & \mbox{for $z\geq L-t$}
    \end{cases}\\
    p(t,z)&=\ 0 \\
    n(t)&=\ 0\\
    a(t,x)&=\  a^0(x+t)\mbox{.}
  \end{split}
\end{displaymath}
For $t>L$ we define inductively $S(t)u=S(L)S(t-L)u$. This
  procedure defines a semigroup of bounded operators in
  $V$ since
  \begin{displaymath}
    \|\psi_1(t,\cdot)\|^2+\|\psi_2(t,\cdot)\|^2+\|a(t,\cdot)\|^2\leq 
    2(1+t^2)^{-\eta}\left(\|\psi_1^0\|+\|\psi_2^0\|+\|a^0\|\right)
  \end{displaymath}
  for $t\leq L$.  The strong continuity of $S$ is a direct
  consequence of the  continuity in the mean in $\Lint^2$.\qed 
\end{pf}
% The operators $S(t)$ have a uniform upper bound 
% \begin{math}
% %  \label{twe:Stup}
%   \|S(t)\|\leq C e^{\gamma t}
% \end{math}
% within finite intervals $[0,T]$.
In order to apply the results of
$C_0$ semigroup theory, we truncate the nonlinearity $g$
smoothly: For any bounded ball $B\subset U$ which is closed w.r.t.\ 
$V$, we choose $g_B:V\to V$ such that $g_B$ is smooth, globally
Lipschitz continuous, and $g_B(u)=g(u)$ for all $u\in B$. This is
possible because the Frechet derivative of $g$ is bounded in $B$ and
the scalar product in $V$ is differentiable with respect to its
arguments.  The corresponding  truncated problem
\begin{equation}
  \label{twe:trunc}
  \frac{d}{dt}u(t)=A u(t)+g_B(u(t))\mbox{,\hspace*{1em}} u(0)=u_0
\end{equation}
has unique global mild and classical solutions (if $u_0\in V$, or
$u_0\in W_{\mathrm{BC}}$, respectively). This implies that the
original problem has the same unique solutions on a sufficiently small
time interval $[0,t_{\mathrm{loc}}]$ for $u_0\in U$.

The following a-priori estimate for the solutions of the truncated
problem \eqref{twe:trunc} extends this local statement to unbounded
time intervals.
\begin{lem}\label{twe:apriori}
  Let $T>0$, $u_0\in U$. If
  $\underline{n}>-\infty$, we suppose $I_k\tau_k>\underline{n}$ for all
  $k=1\ldots m$.  Then, there exists a closed
  bounded ball $B$ such that $B\subset U$ and the solution $u(t)$ of
  the $B$-truncated problem \eqref{twe:trunc} starting at $u_0$ stays
  in $B$ for all $t\in[0,T]$.
\end{lem}
\begin{pf}
  First, let $u_0=(\psi^0, p^0, n^0, a^0)\in D(A)=W_{\mathrm{BC}}\cap U$.  
  
  \emph{Preliminary consideration}\\
  Let $n_*\in(\underline{n},n_k^0)$ be such that
  $G_k(n_*)-\rho_k(n_*)<0$ for all $k=1\ldots m$.  Let $t_1>0$ be such
  that the solution $u(t)=(\psi(t),p(t),n(t),a(t))$ of the
  non-truncated problem \eqref{twe:evolsys} exists in $[0,t_1]$, and
  $n_k(t)\geq n_*$ for all $k=1\ldots m$ and $t\in[0,t_1]$. We define
  the function
  \begin{equation*}
    D(t):=\frac{P}{2}\|\psi(t)\|^2+\sum_{k=1}^m l_k (n_k(t)-n_*)\mbox{.}
  \end{equation*}
  Because of the structure of the nonlinearity $g$, which is linear in
  $\psi$ in its first component, $u(t)$ is classical in $[0,t_1]$.  Hence, $D(t)$ is
  differentiable and the differential equations \eqref{twe:twe} and
  \eqref{twe:carrier} imply
  \begin{align}
    \label{twe:hdiffineq}
    \frac{d}{dt}D(t) &\leq J +\sup_{z\in\C}\{|r_0z+a^0(t)|^2-|z|^2\} - 
    \sum_{k=1}^m \left[\frac{l_k}{\tau_k}
      n_k+ P \re d_k  \int\limits_{S_k}|\psi(z)|^2\,dz\right] \nonumber\\
    &\leq J+\frac{|a^0(t)|^2}{1-|r_0|^2} - \tilde\tau^{-1} n_* -
    \gamma D(t)\mbox{,}       
  \end{align}
  where   (keeping in mind that $\re d_k>0$)
  \begin{displaymath}
  J:=\sum_{k=1}^m l_k I_k\mbox{,\quad}
  \tilde\tau^{-1}:=\sum_{k=1}^m l_k \tau_k^{-1}\mbox{,\quad}
   \gamma:=\min\left\{\tau_k^{-1},\frac{\re d_k}{2}: k=1\ldots
  m\right\}>0\mbox{.} 
  \end{displaymath}
  Consequently, 
  \begin{align}
    \label{twe:hineq}
      D(t) &\leq D(0)+ J t-\tilde\tau^{-1}t\,n_*
      +\frac{1}{1-|r_0|^2}\int_0^t |a^0(s)|^2\,ds\nonumber\\
      &\leq D(0)+ J T+\tilde\tau^{-1}T\,|n_*|+
      \frac{(1+T^2)^{-\eta}}{1-|r_0|^2}\|a^0\|^2\nonumber\\
      &\leq \left(\frac{P}{2}\|\psi^0\|^2+\sum_{k=1}^m l_k n_k^0 + J T +
      \frac{(1+T^2)^{-\eta}}{1-|r_0|^2}\|a^0\|^2\right)
    + \left(L+\tilde\tau^{-1}T\right)\, |n_*| \nonumber\\
    &\leq M+\xi\, |n_*|
  \end{align}
% $h(t)\leq\max\{h(0),\gamma^{-1}J-
%   \gamma^{-1}\tilde\tau^{-1}n_*\}$ for all $t\in[0,t_1]$.  Since
%   $h(0)=\frac{P}{2}\|\psi^0\|^2+ \sum_{k=1}^m l_k n_k^0 - L
%   n_*$, we obtain the estimate
%   \begin{equation}
%     0\leq h(t)\leq M-\xi \cdot n_*
%     \label{twe:hineq}   
%   \end{equation}
  for all $t\in[0,t_1]$ where $M$ and $\xi$ do not depend on $n_*$.
  The inequality \eqref{twe:hineq} remains valid even if $u_0\in
  (V\setminus W_\mathrm{BC})\cap U $ (i.e., $u(t)$ is not classical
  but mild) as both sides of \eqref{twe:hineq} depend only on the
  $V$-norm of $u$ but not on its $W_\mathrm{BC}$-norm. Since
  $n_k(t)\geq n_*$ in $[0,t_1]$ for all $k=1\ldots m$, the estimate
  \eqref{twe:hineq} for $D(t)$ and the differential equation
  \eqref{twe:pol} for $p$ imply bounds for $\psi$, $p$ and $n$ in
  $[0,t_1]$:
  \begin{eqnarray}
    \|\psi(t)\|^2&\leq& S(n_*)^2:=2P^{-1} (M+\xi\,|n_*|) \nonumber\\ 
    \|p(t)\|&\leq& \|p^0\|+S(n_*) \label{twe:bounds}\\ 
    n_k&\in&
    \left[n_*,n_*+ (2 l_k)^{-1} P S(n_*)^2\right]\mbox{.}\nonumber
  \end{eqnarray}
  Hence, $f_k(n_*,(\psi(t),p(t)))$ is greater than 
  \begin{equation}
    \label{twe:nlowineq}
    %f_k(n_*,(\psi(t_1),p(t_1))) > 
    I_k-\frac{n_*}{\tau_k}-
    \frac{P}{l_k}\max_{\Theta\in\R} 
    \left[ (G_k(n_*)-\rho_k(n_*))\Theta^2
    +|\rho_k(n_*)|(|p^0\|+S(n_*))\Theta\right] 
  \end{equation}
  for all $k=1\ldots m$ and $t\in[0,t_1]$. 

  \emph{Construction of $B$}\\ Since
  $G_k(\nu)\to_{\nu\searrow\underline{n}}-\infty$ and $\rho_k(\nu)$ bounded
  for $\nu\searrow\underline{n}$, we can find a $n_*$
  such that the expression \eqref{twe:nlowineq} is greater than $0$
  for all $k=1\ldots m$.  We choose $B$ such that
  $u=(\psi,p,n,a)\in B$ if $\psi$, $p$ and $n$ satisfy
  \eqref{twe:bounds} for this $n_*$ and $a=a^0(t+\cdot)$ for
  $t\in[0,T]$.

  \emph{Indirect proof of invariance of $B$}\\ Assume that the solution
  $v(t)=(\psi(t),p(t),n(t),a(t))$ of the $B$-truncated problem leaves
  $B$. The preliminary consideration and the construction of $B$ imply
  that there exists a $t_1$ such that $u(t)$ exists in $[0,t_1]$, and,
  for one $k\in\{1\ldots m\}$, $n_k(t_1)=n_*$ and $n_k(t)>n_*$ for all
  $t\in[0,t_1]$.  Consequently, $\dot
  n_k(t_1)=f_k(n_k(t_1),(\psi(t_1),p(t_1)))\leq0$.  However, this
  is in contradiction of the construction of $n_*$ such that
  \eqref{twe:nlowineq} is greater than $0$.
%  Now, we prove indirectly
%  that $u(t)$ exists for all $t\geq 0$, and $n_k(t)\in(n_\mathrm{low},
%  n_\mathrm{up}$ for all $k=1\ldots m$ and $t\geq 0$.  Then,
%  $n_k(t)\neq n_\mathrm{low}$ for all $t\geq 0$, and $k=1\ldots m$
%  since otherwise for all satisfy property (P).  Therefore, we can
%  choose the ball $B$ such that the bounds \eqref{twe:bounds} are met
%  by all $w\in B$.
\qquad
\qed\end{pf}

Let us define the map $S(t;u_0):=u(t)$ where $u(t)$ is the mild
solution of the evolution equation \eqref{twe:evolsys} with initial
value $u(0)=u_0$.  The a-priori estimate of Lemma~\ref{twe:apriori}
and standard $C_0$ semigroup theory imply:
\begin{thm}[semiflow]
  \label{twe:existence}\ \\ 
   The map $S$ defines a semiflow mapping $[0,\infty)\times U$ into
  $U$. The map $(t,u_0)\to S(t;u_0)$  is smooth with
  respect to $u_0$ and strongly continuous with respect to $t$.
  If $u_0\in W_{\mathrm{BC}}$ then $S(t;u_0)$ is a classical solution of
  \eqref{twe:evolsys} for all $t\geq0$.
\end{thm}
The smooth dependence of the solution on all parameters within a
boun\-ded parameter region is also a direct consequence of the $C_0$
semigroup theory. The restrictions imposed on the parameters in
Section~\ref{sec:twe} and Lemma~\ref{twe:apriori} have to be satisfied
uniformly in the parameter range under consideration in order to
obtain a uniform a-priori estimate.

Theorem~\ref{twe:existence} still permits for $S(t;u_0)$ to grow
to $\infty$ for $t\to\infty$. The following corollary observes that
this is not the case if the component $a^0$ of $u_0$ is globally
bounded, that is, $a^0\in L^\infty$. This is of practical importance
because $a^0$, representing the optical injection into the laser, is
always bounded. A physically sensible model should provide globally
bounded solutions in this case.
\begin{cor}[global boundedness]
  \label{twe:boundedness}
  Let $u_0=(\psi^0,p^0,n^0,a^0)\in U$ where, in addition,
  $\|a^0\|_\infty<\infty$.  Then there exists a constant $C$ such that
  $\|S(t;u_0)\|_V\leq C$ for all $t\geq0$.
\end{cor}
\begin{pf}
  It is sufficient to prove that the constants $M$ and $\xi$ in the
  estimate \eqref{twe:hineq} for $D(t)$ do not depend on $T$ if
  $\|a^0\|_\infty<\infty$. The estimate \eqref{twe:hdiffineq} for $\frac{d}{dt} D(t)$ implies
  \begin{align}
    \label{twe:hineqb}
    D(t)&\leq \max\left\{D(0),\frac{1}{\gamma}\left(J+\frac{\|a^0\|_\infty}{1-|r_0|^2}-
      \frac{n_*}{\tilde\tau}\right)\right\}\nonumber\\
  &\leq \left(\frac{P}{2}\|\psi^0\|^2+\sum_{k=1}^m l_k n_k^0 +L |n_*|\right)
  +\frac{1}{\gamma}\left(J+\frac{\|a^0\|_\infty}{1-|r_0|^2}+
    \frac{|n_*|}{\tilde\tau}\right)\nonumber\\
  &\leq \left(\frac{P}{2}\|\psi^0\|^2+\sum_{k=1}^m l_k n_k^0+
    \frac{1}{\gamma}\left[J+\frac{\|a^0\|_\infty}{1-|r_0|^2}\right]\right)+
  \left(L+\frac{1}{\gamma\tilde\tau}\right)\,|n_*|\nonumber\\
  &\leq M+\xi\,|n_*|
  \end{align}
  where now $M$ and $\xi$ do not depend on $T$. Hence, the bounds
  \eqref{twe:bounds} can now be derived from \eqref{twe:hineqb} in the
  same way as in the proof of Lemma~\ref{twe:apriori} using the
  $T$-independent bounds $M$ and $\xi$. Consequently, we can choose
  $n_*$ independent of $T$ and, hence, the ball $B$ does not depend on
  $T$ (see proof of Lemma~\ref{twe:apriori}).
\qed\end{pf}

\section{Asymptotic behavior of the linear part --- 
  spectral gap for $H(n)$}
\label{sec:spec}

We restrict ourselves to the autonomous system
\eqref{twe:twe}--\eqref{twe:carrier} in the following. The boundary
conditions are
\begin{equation}
  \label{spec:bcond}
   \psi_1(t,0)=r_0\psi_2(t,0)\mbox{,\quad }
  \psi_2(t,L)=r_L\psi_1(t,L)
\end{equation}
in the autonomous case. 

As mentioned in Section~\ref{sec:nt}, the long-time behavior of the overall
system at $\epsilon=0$ in \eqref{spec:carrier} (i.e., $\dot n_k=0$ for
$k=1\ldots m$) is determined by the behavior of the linear space-dependent
subsystem \eqref{spec:linsub}, that is, the spectral properties of the
operator $H(n)$ defined by \eqref{spec:h}.  We treat $n$ as a parameter in
this section, aiming to find so-called \emph{critical carrier densities} $n$ that
give a spectral splitting for the strongly continuous semigroup $T(n;t)$
generated by $H(n)$.
\begin{defn}[Critical density]\label{spec:critn}
  Let us denote a carrier density $n$ as \emph{critical} if there exists
   a splitting of $X$ into two $H(n)$-invariant
  subspaces $X_c(n)\oplus X_s(n)=X$ such that $X_c(n)$ is
  finite-dimensional, $H(n)\vert_{X_c(n)}$ has all eigenvalues on the
  imaginary axis, and
  \begin{displaymath}
    \|T(n;t)\vert_{X_s(n)}\|\leq M e^{-\xi t}\mbox{\qquad for some $M>0$, $\xi>0$
    and all $t\geq0$.}
  \end{displaymath}
\end{defn}
The general result of \cite{NRL86} implies that $T(n;t)$ is a compact
perturbation of the semigroup $T_0(n;t)$ generated by the diagonal part
$H_0(n)$ of $H(n)$ where
\begin{displaymath}
  H_0=
  \begin{bmatrix}
    \begin{bmatrix}
      -\partial_z+\beta(n) &0\\ 0& \partial_z+\beta(n)
    \end{bmatrix} & 0\\[0.5em] 0 & i\Omega_r(n)-\Gamma(n) 
  \end{bmatrix}
\end{displaymath}
is defined in $Y\subset X$ and $\beta(n)=(1+i\alpha_H)G(n)-d-\rho(n)$.
The growth rate of $T_0$ is governed by 
\begin{equation}\label{spec:t0}
  \|T_0(n;t)\|\leq M \exp(R_\infty(n)\, t)
\end{equation}
with $R_\infty(n)=\ \max\left\{R_\psi(n),R_p(n)\right\}$ where 
\begin{displaymath}
  R_\psi(n)= \frac{1}{L}\left[
      \sum_{k=1}^ml_k\re\beta_k(n_k)+\frac{1}{2}\log|r_0r_L|\right]\mbox{,\quad}
    R_p(n)=\ -\min_{k=1,\ldots,m}\{\Gamma_k(n_k)\}\mbox{.}
\end{displaymath}
The quantity $R_\psi$ gives the growth rate for the first two
components of $T_0$ (corresponding to $\psi$), whereas $R_p$ limits the growth of
the components corresponding to $p$.  In particular, the first two
components of $T_0$ decay to zero after time $2L$ if $r_0r_L=0$.
Moreover, $R_\psi(n)$ tends to $-\infty$ for $n\to\underline{n}$ also
for $r_0r_L\neq0$ because $G_k(n_k)\to-\infty$ if
$n_k\searrow\underline{n}$ for all $k=1,\ldots,m$, and $R_p(n)<-1$ for
all $n$ (see Section~\ref{sec:twe} for fundamental assumptions on the
parameters). The estimate \eqref{spec:t0} is sharp for
$T_0$. The quantity $\exp(tR_\infty(n))$ determines the radius of the
essential spectrum of $T(n;t)$

The following lemma establishes that the presence of a nonzero
$\kappa_k$ or a positive $\rho_k(\cdot)$ makes it possible to find a critical
carrier density $n$.
\begin{lem}[Existence of critical carrier density with spectral gap]
  \label{spec:lemgap}
  \ \\  Assume that one of the following two conditions is satisfied.
  \begin{enumerate}
  \item\label{kappa} The coupling $\kappa_k$ is nonzero for at least
    one $k\in\{1,\ldots,m\}$, and $m\geq2$.
  \item\label{rho} The reflectivities satisfy $r_0r_L\neq0$,
    $\kappa_k=0$ for all $k\in\{1,\ldots,m\}$, and $\rho_k(\nu)>0$ for
    all $\nu>1$ for at least one $k\in\{1,\ldots,m\}$.
  \end{enumerate}
  Then there exists a critical carrier density $n$
  in the sense of Definition~\ref{spec:critn}.  
\end{lem}
\begin{pf}\emph{(i) Preliminary observation}\\
  Since $T(n;t)$ is a compact perturbation of $T_0(n;t)$ for all
  $t\geq0$, $T(n;t)$ cannot decay with a rate faster than
  $-R_\infty(n)$. Moreover, $H(n)$ can have at most finitely many
  eigenvalues with real part bigger than $R_\infty(n)+\delta$ for any
  given $\delta>0$. We will first derive the characteristic function
  $h(n;\cdot)$, the roots of which are the eigenvalues of $H(n)$.
  The statement of the lemma can then be proved by finding a $n$ such
  that $R_\psi(n)<0$, and $h(n;\cdot)$ has no roots with positive real
  part and at least one root on the imaginary axis.
  
  Let $\lambda$ be an eigenvalue of $H(n)$ and $(\psi,p)$ its
  eigenvector. For a general $n\in(\underline{n},\infty)^m$ define the
  functions $\chi_k:(\underline{n},\infty)\times\C \to \C$ by
  \begin{displaymath}
    \chi_k(n_k;\lambda)=\frac{\rho_k(n_k)\Gamma_k(n_k)}{\lambda-
      i\Omega_{r,k}(n_k)+\Gamma_k(n_k)}
  \end{displaymath}
  having poles at $i\Omega_{r,k}(n_k)-\Gamma_k(n_k)$ ($k=1,\ldots,m$).
  The component $\psi$ of the eigenvector satisfies the boundary value
  problem 
  \begin{equation}
    \label{spec:resolvbvp}
    \begin{gathered}
      \begin{bmatrix}
        -\partial_z +\beta(n)+\chi(n;\lambda)-\lambda & -i\kappa\\
        -i\kappa & \partial_z +\beta(n)+\chi(n;\lambda)-\lambda
      \end{bmatrix}
      \begin{bmatrix}
        \psi_1\\ \psi_2
      \end{bmatrix}
      =0\\[0.5em]
      \mbox{ with b.c.\quad}\psi_1(0)=r_0\psi_2(0)\mbox{,\quad }
      \psi_2(L)=r_L\psi_1(L)\qquad
    \end{gathered}
  \end{equation}
  on the interval $[0,L]$. The coefficients $\beta(n)$, $\kappa$, and
  $\chi(n;\lambda)$ in \eqref{spec:resolvbvp} are piecewise constant
  in $z$. That is, for $z\in S_k$,
  $[\beta(n)\psi_1](z)=\beta_k(n_k)\psi_1(z)$,
  $[\kappa\psi_1](z)=\kappa_k\psi_1(z)$, and
  $[\chi(n;\lambda)\psi_1](z)=\chi_k(n_k;\lambda)\psi_1(z)$ (and
  likewise for $\psi_2$).  Denote the overall transfer matrix of
  \eqref{spec:resolvbvp} by $T(z_1,z_2,n;\lambda)$ for
  $z_1,z_2\in[0,L]$. Then $\psi(z)$ is a multiple of $T(z,0,n;\lambda)
  \left(\begin{smallmatrix} r_0\\ 1
    \end{smallmatrix}\right)$, and 
  $p=\Gamma(n)\psi/(\lambda-i\Omega_r(n)+\Gamma(n))$.  Thus, the
  eigenvalue $\lambda$ is geometrically simple if
  $\lambda\neq i\Omega_{r,k}(n_k)-\Gamma_k(n_k)$ ($k=1,\ldots,m$).
  
  The fact that all coefficients of \eqref{spec:resolvbvp} are
  piecewise constant in $z$ makes it possible to find the transfer
  matrix $T(L,0,n;\lambda)$ corresponding to \eqref{spec:resolvbvp}
  analytically; see \cite{B94}, \cite{RSS97}. Within each subinterval
  $S_k$ the transfer matrix is given by
  \begin{equation}
    \label{spec:transfer}
    T_k(z,n_k;\lambda)=\frac{e^{-\gamma_k z}}{2\gamma_k}
    \begin{bmatrix}
      \gamma_k+\mu_k+e^{2\gamma_k z}[\gamma_k-\mu_k] &
      i \kappa_k \left[1-e^{2\gamma_k z}\right] \\
      -i \kappa_k \left[1-e^{2\gamma_k z}\right] &
      \gamma_k-\mu_k+e^{2\gamma_k z}[\gamma_k+\mu_k]
    \end{bmatrix}
  \end{equation}
  where
  $\mu_k=\mu_k(n_k;\lambda)=\lambda-\chi_k(n_k;\lambda)-\beta_k(n_k)$
  and $\gamma_k=\gamma_k(n_k;\lambda)=\sqrt{\mu_k^2+\kappa_k^2}$. The
  right-hand-side of \eqref{spec:transfer} does not depend on the
  branch of the square root in $\gamma_k$ since the expression is even
  with respect to $\gamma_k$.  The function
  \begin{equation}\label{spec:hchar}
    h(n;\cdot)=
    \begin{bmatrix}
      r_L, &-1
    \end{bmatrix}T(L,0,n;\cdot)
    \begin{bmatrix}
      r_0\\ 1
    \end{bmatrix}=
    \begin{bmatrix}
      r_L &-1
    \end{bmatrix}
    \prod_{k=m}^1 T_k(l_k,n_k;\cdot)
    \begin{bmatrix}
      r_0\\ 1
    \end{bmatrix}
  \end{equation}
  defined in $\C\setminus\{i\Omega_{r,k}(n_k)-\Gamma_k(n_k):
  k=1,\ldots,m\}$ is the characteristic function of $H(n)$. Its roots are
  the eigenvalues of $H(n)$. The algebraic multiplicity of any eigenvalue
  equals the multiplicity of the corresponding root.
  
  \emph{(ii) Coarse upper bound for eigenvalues}\\
  We prove in the first step that all eigenvalues $\lambda$ of
  $H(n)$ satisfy
  \begin{equation}\label{spec:relambda}
    \re\lambda<\Lambda_u(n):=\max_{k=1\ldots m}
    \left\{\re \beta_k(n_k)+2\rho_k(n_k),-\frac{\Gamma_k(n_k)}{2}\right\}\mbox{.} 
  \end{equation}
  Partial integration of the eigenvalue equation
  \eqref{spec:resolvbvp} and its complex conjugate equation yields
  \begin{displaymath}
    \re\lambda\leq \max_{k=1\ldots m}\left\{\re\chi_k(n_k;\lambda)
      +\re \beta_k(n_k)\right\}
  \end{displaymath}
  for any eigenvalue $\lambda$ of $H(n)$.  If
  $\re\lambda\geq-\Gamma_k(n_k)/2$ we get
  $\re\chi_k(n_k;\lambda)\leq|\chi_k(n_k;\lambda)|\leq 2\rho_k(n_k)$,
  thus, \eqref{spec:relambda} holds for eigenvalues of $H(n)$.
  
  We observe that $\Lambda_u(n)<0$ if $n_k$ are close to
  $\underline{n}$ for all $k\in\{1,\ldots,m\}$, as
  $\rho_k(n_k)/G_k(n_k)\to0$ for $n_k\searrow\underline{n}$ and
  $\Gamma_k(n_k)>1$.
  
  \emph{(iii) Proof for a single section, $m=1$, with $r_0=0$}\\
  We first prove the statement of the lemma for the case $m=1$,
  $r_0=0$. Then we treat the multi-section case ($m>1$) as a
  perturbation of the single section by choosing $n_k$ sufficiently
  close to $\underline{n}$ for all other sections. For brevity, we
  drop the index $1$ in this paragraph since all quantities are
  one-dimensional.  There are at most finitely many eigenvalues with
  real parts bigger that $R_\infty(n)+1/2$, which is uniformly
  negative for all $n$ due to $r_0=0$ ($r_0=0$ implies $R_\psi(n)=-\infty$,
  and, hence, $R_\infty(n)=R_p(n)$). Thus, it is sufficient to find an
  $n$ such that $h(n;\cdot)$ has a root with positive real part.  Then the
  $n$ satisfying the statement of the lemma must exist and be less or
  equal to $n$ due to observation (ii).

  For $m=1$, $\kappa\neq0$, $r_0=0$ the characteristic function $h$
  admits the form
  \begin{multline}\label{spec:htilde}
    \tilde
    h(\mu)=r_L\frac{i\kappa}{2\sqrt{\mu^2+\kappa^2}}
    \left[\exp\left[-2l\sqrt{\mu^2+\kappa^2}\right]-1\right]-\\ 
    -\frac{\sqrt{\mu^2+\kappa^2}+\mu}{2\sqrt{\mu^2+\kappa^2}}
    -\frac{\sqrt{\mu^2+\kappa^2}-\mu}{2\sqrt{\mu^2+\kappa^2}}
    \exp\left[-2l\sqrt{\mu^2+\kappa^2}\right]
  \end{multline}
  after multiplication with the nonzero factor
  $\exp\left[-l\sqrt{\mu^2+\kappa^2}\right]$ where
  $\mu=\mu(n;\lambda)=\lambda+d-[1+i\alpha_{H}]G(n)+\rho(n)-\chi(n;\lambda)$.
  Since $G$ is monotone increasing to infinity and dominating $\rho$
  for large $n$, the equation $\tilde \mu=\mu(n;i\omega+\delta)$ can
  be solved for $(n,\omega)\in(1,\infty)\times\R^+$ for any given
  $\delta>0$ if $-\re\tilde\mu$, $\im\tilde\mu$ and $-\im\tilde
  \mu/\re\tilde\mu$ are sufficiently large.  This implies that, if $\tilde
  h(\cdot)$ has roots $\tilde\mu$ with a sufficiently large
  $-\re\tilde\mu$, $\im\tilde\mu$ and $-\im\tilde\mu/\re\tilde\mu$, we
  can find a $n$ such that $\tilde h(\mu(n;\cdot))$ has a root
  $\lambda$ in the positive half-plane.
  
  The function $\tilde h$ has at least one sequence of roots
  $\mu_j$ with $|\mu_j|\to_{j\to\infty}\infty$. More precisely, the
  roots satisfy $\mu_j-\mu_j^0\to_{j\to\infty}0$ where
  $\mu_j^0=W_j(i\kappa l/2)/l$ for $r_L=0$, and $\mu_j^0=W_j(ir_L\kappa
  l)/(2l)$ for $r_L\neq0$. The function $W_j$ denotes the $j$-th
  complex sheet of the Lambert W function (defined by $W_j(x)\exp
  W_j(x)=x$) in the upper complex half-plane. The convergence (and
  existence) of the roots $\mu_j$ of $\tilde h$ is a direct
  consequence of the fact that the sheets $W_j$ correspond to
  uniformly simple roots of an analytic function and that the sequence
  of functions
  $$\zeta\mapsto\tilde
  h(\mu_j^0+\zeta)-\left[\frac{\kappa\exp(-2l(\mu_j^0+\zeta))}{2(\mu_j^0+\zeta)}
    \left(ir_L-\frac{\kappa}{2(\mu_j^0+\zeta)}\right)-1\right]$$
  converges to $0$ uniformly on a sufficiently small ball around $0$.
  
  The complex numbers $\mu_j^0$ (and, thus, $\mu_j$) satisfy
  $\re\mu_j^0\to-\infty$, $\im\mu_j^0\to+\infty$ and
  $-\im\mu_j^0/\re\mu_j^0\to+\infty$. Hence, we can find $j$,  $n$
  and  $\lambda$ such that $\mu_j=\mu(n,\lambda)$, $\re\lambda>0$ and
  $\im\lambda>0$. This $n$ and $\lambda$ satisfy $h(n;\lambda)=0$.
  Consequently, $\lambda$ is an eigenvalue of $H(n)$ with positive
  real and imaginary part.

  \emph{(iv) Proof for case (\ref{kappa})}\\
  Assume without loss of generality that the section $S_{\tilde k}$
  with nonzero $\kappa$ has an index $\tilde k\geq2$.  

  Let $k\neq \tilde k$ be
  in $\{1,\ldots,m\}$. We observe that
  $\exp[-l_k\gamma_k(n_k;\lambda)]\,T_k(l_k,n_k;\lambda)$ behaves
  asymptotically for $n_k\searrow \underline{n}$ like
  \begin{displaymath}
    \exp[-l_k\gamma_k(n_k;\lambda)]\,T_k(l_k,n_k;\lambda) =
    \begin{bmatrix}
      0&0\\[-1ex] 0&1
    \end{bmatrix}+R_k(n_k;\lambda)\mbox{.}
  \end{displaymath}
  The remainder $R_k(n_k;\lambda)$ tends to zero for
  $n_k\to\underline{n}$ uniformly with respect to $\lambda$ in any strip
  of finite width around $i\R^+$, as well as its derivative with respect
  to $\lambda$. Consequently,  the function
  $$\hat h(n;\cdot)=\exp(-\sum_{j=1}^m\gamma_j(n_j;\cdot)l_j)\,
  h(n;\cdot)\mbox{,}$$
  a nonzero multiple of $h(n;\cdot)$, converges
  to $\tilde h(\mu_{\tilde k}(n_{\tilde k};\cdot))$ (where $r_L$ has
  to be replaced by $0$ if $\tilde k<m$ in \eqref{spec:htilde}) when
  all $n_k$ with $k\neq \tilde k$ tend to $\underline{n}$. This
  convergence is uniform in any strip of finite width around $i\R^+$.
  
  Observation (iii) has shown that there exists a $n_{\tilde k}$ such
  that $\tilde h((\mu_{\tilde k}(n_{\tilde k};\cdot))$ has a root
  $\lambda$ with positive real and imaginary part. This root stays in
  the positive half-plane under the perturbation toward $\hat
  h(n;\cdot)$ if all $n_k$ with $k\neq \tilde k$ are sufficiently
  close to $\underline{n}$.  Furthermore, $R_\psi(\nu)<-1$ for all
  $\nu\in(\underline{n},\infty)^m$ with $\nu_{\tilde k}\leq n_{\tilde
    k}$ and $\nu_k$ sufficiently close to $\underline{n}$ for all
  $k\neq \tilde k$.  If $\nu_{\tilde k}$ is also close to
  $\underline{n}$ then the upper bound for eigenvalues
  $\Lambda_u(\nu)<0$ due to observation (ii). Consequently, since
  eigenvalues of $H(\nu)$ depend continuously on $\nu$, there must
  exist a $\nu$ satisfying the statement of the lemma with
  $\nu_{\tilde k}\in(\underline{n}, n_{\tilde k})$ and $\nu_k$
  sufficiently close to $\underline{n}$ for all $k\neq \tilde k$.
  
  \emph{(v) Proof for case (\ref{rho})}\\
  If all components of $n$ are equal to $1$ then $R_\psi(n)<0$.  Due
  to $r_0r_L\neq0$, $R_\psi(n)$ tends to infinity if one component of
  $n$ tends to infinity while the others remain fixed. It is
  sufficient to show that $h(n;\lambda)$ has at least one root with
  real part bigger than $R_\psi(n)$ if $n$ is such that
  $R_\psi(n)>-1$. (The radius of the essential spectrum of $T(n;t)$ is
  $\exp(t\max\{R_\psi(n),R_p(n)\})$ and $R_p(n)<-1$ for all $n$.)

  Since $\kappa_k=0$ for all $k=1,\ldots,m$, the characteristic
  equation $h(n;\lambda)=0$ simplifies to 
  \begin{displaymath}
    r_0r_L=\exp(2\sum_{k=1}^ml_k\mu_k(n_k;\lambda))
  \end{displaymath}
  A complex number $\lambda$ is a solution of this equation if (and
  only if) there exists a $j\in\Z$ such that 
  \begin{align}\label{spec:fp}
      \lambda&=\  \frac{1}{L}\left[
        \sum_{k=1}^ml_k\mu_k(n_k;\lambda)+\frac{1}{2}\log(r_0r_L)+j\pi i\right]\nonumber\\
      &=\ \frac{1}{L}\left[
        \sum_{k=1}^ml_k\left[\beta_k(n_k)+\chi_k(n_k;\lambda)\right]+\frac{1}{2}\log(r_0r_L)+j\pi
        i\right]\nonumber\\
      &=\ \lambda_0(n)+\frac{j\pi i}{L}+\sum_{k=1}^m\frac{l_k}{L}\chi_k(n_k;\lambda)
  \end{align}
  where
  $\lambda_0(n)=\frac{1}{L}\left[\frac{1}{2}\log(r_0r_L)+\sum_{k=1}^ml_k\beta_k(n_k)\right]$
  (note that $\re\lambda_0(n)=R_\psi(n)$). For large $j$, equation
  \eqref{spec:fp} has exactly one solution $\lambda_j(n)$ close to
  $\lambda_0(n)+\frac{j\pi i}{L}$ because $\chi_k(n_k;\lambda)\to0$
  for $\im\lambda\to\infty$ for all $k\in\{1,\ldots,m\}$. Moreover,
  \eqref{spec:fp} can be used to obtain this solution by fixed point
  iteration starting from $\lambda_0(n)+\frac{j\pi i}{L}$ since
  $\partial_\lambda\chi_k(n_k;\lambda)\to0$ for $\im\lambda\to\infty$
  for all $k\in\{1,\ldots,m\}$. The real part of $\chi_k(n_k;\lambda)$
  is
  \begin{displaymath}
    \re\chi_k(n_k;\lambda)=\frac{\rho_k(n_k)\Gamma_k(n_k)
      [\re\lambda+\Gamma_k(n_k)]}{[\re\lambda+\Gamma_k(n_k)]^2
      +[\im\lambda-i\Omega_{r,k}]^2}\mbox{,}    
  \end{displaymath}
  which is non-negative for all $k\in\{1,\ldots,m\}$ and strictly
  positive for at least one $\tilde k$ if $\re\lambda>-1$,
  $n_k\geq1$ for all $k\in\{1,\ldots,m\}$, and $n_{\tilde k}>1$ (due
  to condition (\ref{rho}): $\rho_{\tilde k}(n_{\tilde k})>1$ if
  $n_{\tilde k}>1$). Hence, for large $j$, the result $\lambda_j(n)$ of
  iterating \eqref{spec:fp} has a real part strictly larger than
  $\re\lambda_0(n)$ if $\re\lambda_0(n)>-1$.  Thus,
  $\re\lambda_j(n)>R_\psi(n)=\re\lambda_0(n)$ if $n$ is such that
  $R_\psi(n)>-1$.  \qed
\end{pf}
\paragraph*{Remarks} The conditions of Lemma~\ref{spec:lemgap} cover all
practically relevant cases except the case of a single section laser with
nonzero $\kappa$ and $r_0r_L\neq 0$. The case $r_0r_L=0$, $m=1$, $\kappa\neq0$
is not mentioned in Lemma~\ref{spec:lemgap} but also proven in part (iii) of
its proof. The conditions of Lemma~\ref{spec:lemgap} also highlight the two
effects causing a spectral gap of the linear wave operator, the nonzero
coupling (nonzero $\kappa$), and the gain dispersion (positive $\rho$).

The set of critical densities $n$ consists piecewise of smooth submanifolds
of $\R^m$.   The question whether there is always point spectrum to the right
of the essential radius $R_\infty(n)$  in the case of
nonzero coupling $\kappa_k$ and $r_0r_L\neq0$ is open.

\section{Existence and properties of the finite-di\-men\-sio\-nal  
  center ma\-ni\-fold}
\label{sec:im}
In this section we construct a low-dimensional attracting invariant
manifold for system \eqref{spec:linsub}, \eqref{spec:carrier} using
the general theorems about the persistence and properties of normally
hyperbolic invariant manifolds in Banach spaces \cite{BLZ98},
\cite{BLZ99}, \cite{BLZ00}. The statements of the theorem and the
proofs rely only on the system's structure
\begin{equation}
  \label{im:struc}
  \begin{split}
  \frac{d}{dt}E&=\ H(n)E\\    
  \frac{d}{dt}\,n&=\ \epsilon F(n,E)
  \end{split}
\end{equation}
the spectral properties of $H(n)$ for fixed critical densities $n$, the
smoothness of the semiflow $S(t;\cdot)$ generated by \eqref{im:struc} with
respect to parameters and initial values, the smallness of $\epsilon$,
and that $\partial_2 F(n,0)=0$.

In particular, we investigate system \eqref{im:struc} in the vicinity of
critical carrier densities as constructed in Lemma~\ref{spec:lemgap}, assuming
throughout this section that one of the conditions of Lemma~\ref{spec:lemgap}
is satisfied. Let $\mathcal{K}\subset\R^m$ be a set of critical carrier
densities with a uniform spectral gap of size greater than $\xi$. Then for any
$\delta\in(0,\xi]$ there exists a simply connected open neighborhood
$U_\delta$ of $\mathcal{K}$ such that
\begin{align*}
  \spec H(n)&=\ \sigma_c(n)\cup\sigma_s(n) &&\mbox{where}\\
  \re\sigma_c(n)&>\ -\delta\mbox{,} &&\\
  \re\sigma_s(n) &<-\xi &&\mbox{for all $n\in \clos U_\delta$.}
\end{align*}
The number of elements of $\sigma_c(n)$ is finite and, hence, constant in
$U_\xi$ if the eigenvalues are counted according to their algebraic
multiplicity. We denote this number by $q$.  Note that $\mathcal{K}$ is
typically a submanifold of dimension $m-q$ in $\R^m$. In the special case
$m=q$ the invariant manifold, constructed in this section, corresponds to a
local center manifold.

There exist spectral projections of $H(n)$, $P_c(n)$ and
$P_s(n)\in\mathcal{L}(X)$, corresponding to this splitting. They are
well defined and unique for all $n\in U_\xi$ and depend smoothly on
$n$. We define the corresponding closed invariant subspaces of $X$ by
$X_c(n)=\im P_c(n)=\ker P_s(n)$ and $X_s(n)=\im P_s(n)=\ker P_c(n)$.
The complex dimension of $X_c(n)$ is $q$. Let $B(n):\C^q\to X$ be a
basis of $X_c(n)$ which depends smoothly on $n$.  $B(\cdot)$ is
well defined in $U_\xi$ because $U_\xi$ is simply connected, has
rectifiable boundary and $H$ has a uniform spectral splitting on
$\clos U_\xi$.  The existence of the basis $B$ and the spectral
projection $P_c$ and their smooth dependence on $n\in U_\xi$ imply
that the maps $\tilde P_c: U_\xi\mapsto \mathcal{L}(X;\C^q)$ and
$\mathcal{ R}:X\times U_\xi \mapsto \C^q\times U_\xi$ defined by
\begin{displaymath}
  \tilde P_c(n)=B(n)^{-1}P_c(n)\mbox{,}\quad \mathcal{
  R}(E,n):=(\tilde P_c(n)E,n)
\end{displaymath}
are well defined and smooth. We also know that the semiflow
$T(n;\cdot)\vert_{X_s(n)}=P_s(n)T(n;\cdot)$ decays with a rate strictly
greater than $\xi$ for all $n\in U_\xi$.
Using these notations, we can state the following
theorem about the existence and properties of invariant manifolds of
the semiflow $S(t;\cdot)$ of system \eqref{im:struc}:
\begin{thm}[Model reduction]
  \label{im:thm}
  \ \\ Let $k>2$ be an integer number, $\delta\in(0,\xi)$ be sufficiently
  small, and $\mathcal{ N}$ be a closed bounded cylinder in $\C^q\times
  U_\delta$. That is, $\mathcal{N}=\{(E_c,n)\in\C^q\times U_\delta:\|E_c\|\leq
  R, n\in\mathcal{N}_n\}$ for a arbitrary $R>0$ and some closed bounded set
  $\mathcal{N}_n\subset U_\delta$ .  Then, there exists an $\epsilon_0>0$ such
  that the following holds.  For all $\epsilon\in[0,\epsilon_0)$, there exists
  a $C^k$ manifold $\mathcal{C}\subset X\times\R^m$ satisfying:
  {\renewcommand{\theenumi}{\roman{enumi}}
    \begin{enumerate}
    \item \emph{(Invariance)} $\mathcal{C}$ is $S(t,\cdot)$-invariant
      relative to $\mathcal{ R}^{-1}\mathcal{N}$. That is, if
      $(E,n)\in\mathcal{ C}$, $t\geq0$, and $S([0,t];(E,n))\subset \mathcal{ R}^{-1}\mathcal{
        N}$, then $S([0,t];(E,n))\subset\mathcal{ C}$.
    \item \emph{(Expansion in $\epsilon$)} $\mathcal{C}$ can be represented as the
      graph of a map which maps 
      \begin{displaymath}
      (E_c,n,\epsilon)\in\mathcal{N}\times[0,\epsilon_0) \to
      ([B(n)+\epsilon\nu(E_c,n,\epsilon)]E_c,n)\in X\times\R^m
      \end{displaymath}
      where $\nu: \mathcal{N}\times[0,\epsilon_0)\to \mathcal{L}(\C^q;X)$
      is $C^{k-2}$ with respect 
      to all arguments. Denote the $X$-component of $\mathcal{C}$ by
      \begin{displaymath}
        E_X(E_c,n,\epsilon)=[B(n)+\epsilon\nu(E_c,n,\epsilon)]E_c\in X\mbox{.}        
      \end{displaymath}
    \item\label{im:attr} \emph{(Exponential attraction/foliation)} Let
      $\Upsilon\subset \mathcal{R}^{-1}\mathcal{N}\subset X\times\R^m$ be a
      bounded set such that $\mathcal{ R}\Upsilon\subset \mathcal{ N}$ has a
      positive distance to the boundary of $\mathcal{ N}$. For any
      $\Delta\in(0,\xi)$ there exist a constant $M$ and a time $t_c\geq0$ with
      the following property: For any $(E,n)\in\Upsilon$ there exists a
      $(E_c,n_c)\in\mathcal{ N}$ such that
%       Let $(E,n)$ be such that
%       $S(t;(E,n))\in\mathcal{N}$ for all $t\geq 0$. Then, there exist
%       $(E_c,n_c)\in\mathcal{B}$, $M>0$ and $t_c\geq 0$ such that
      \begin{equation}
        \label{im:attreq}
        \|S(t+t_c;(E,n))-S(t;(E_X(E_c,n_c,\epsilon),n_c))\|\leq M e^{-\Delta
        t}
      \end{equation}
      for all $t\geq 0$ with $S([0,t+t_c];(E,n))\subset\Upsilon$.
    \item\label{im:flow} \emph{(Flow)} 
      The values
      $\nu(E_c,n,\epsilon)E_c$ are in $D(H(n))$ and their $P_c(n)$-component
      is $0$ for all $(E_c,n,\epsilon)\in
      \mathcal{N}\times[0,\epsilon{}_0)$. 
      The flow on $\mathcal{C}\cap
      \mathcal{ R}^{-1}\mathcal{N}$ is differentiable with respect to $t$ and governed
      by the following system of ODEs:
      \begin{equation}
        \label{im:flowc}
        \begin{split}
          \frac{d}{dt} E_c &=\ \left[H_c(n)+\epsilon
          a_1(E_c,n,\epsilon) +\epsilon^2a_2(E_c,n,\epsilon)
          \nu(E_c,n,\epsilon)\right] E_c\\ 
          \frac{d}{dt} n_{\phantom{k}} &=\ \epsilon F(n,E_X(E_c,n_c,\epsilon))
        \end{split}
      \end{equation}
      where
      \begin{displaymath}
        \begin{split}
          H_c(n) &=\ \tilde P_c(n)H(n)B(n)\\
          a_1(E_c,n,\epsilon) &=\ 
          -\tilde P_c(n)\partial_nB(n)F(n,E_X(E_c,n_c,\epsilon))\\
          a_2(E_c,n,\epsilon) &=\ 
          B(n)^{-1}\partial_nP_c(n)F(n,E_X(E_c,n_c,\epsilon))\,(\id-P_c(n))\mbox{.}
        \end{split}
      \end{displaymath}
      System \eqref{im:flowc} is symmetric with respect to rotation $E_c\to
      E_ce^{i\phi}$ and $\nu$ satisfies the relation
      $\nu(e^{i\phi}E_c,n,\epsilon)=\nu(E_c,n,\epsilon)$ for all
      $\phi\in[0,2\pi)$. 
    \end{enumerate}
    }
\end{thm}

\emph{Remark:} This theorem is an application of the general theory
about persistence of normally hyperbolic invariant manifolds of
semiflows under $C^1$ small perturbations \cite{BLZ98}, \cite{BLZ99},
\cite{BLZ00}. In this case, we find the unperturbed invariant
manifold, which is even finite-dimensional and exponentially stable,
for $\epsilon=0$. The proof describes in detail
the appropriate cut-off modification of the system outside of the
region of interest to make the unperturbed invariant manifold compact.
Then it shows how the results of the previous sections guarantee the
$C^1$-smallness of the perturbation and the normal hyperbolicity.

A model reduction for systems of ODEs with the structure
\eqref{int:geneq} has been presented already by \cite{T01} using
Fenichel's Theorem \cite{F79}.
\begin{pf}
  \emph{\ \\ Cut-off modification of system \eqref{im:struc} outside
    the region    of interest}\\
  The projections $P_c(n)$, $\tilde P_c$ and $P_s(n)$, and the basis
  $B(n)$ are defined only for $n\in U_\xi$. First, we define the
  appropriate $\delta>0$ and extend the definitions of $B$, $P_c$,
  $\tilde P_c$ and $P_s$ to the whole $\R^m$. Let
  $r_1:[0,\infty]\mapsto[0,1]$ be a smooth Lipschitz continuous
  function that satisfies $r_1(x)=0$ for $x<R$ and $r_1(x)=1$ for
  $x>R+1$.  The $\C^q$ map $E_c\mapsto r_1(\|E_c\|)E_c$ is smooth with
  respect to $E_c$ (if $\C^q$ is identified with $\R^{2q}$) and
  Lipschitz continuous with a Lipschitz constant $L$.  We choose the
  $\delta>0$ for the theorem such that
  \begin{equation}\label{im:deltadef}
    \delta<\delta_0:=\frac{\xi}{k(1+L)}\mbox{.}
  \end{equation}
  For fixed $n\in U_{\delta}$ denote by $T_{c,n}^1:\R\times\C^q\mapsto\C^q$
  the nonlinear flow generated by $\dot
  E_c=[H_c(n)+\delta_0r_1(\|E_c\|)]E_c$. The
  definition~\eqref{im:deltadef} of $\delta$ implies:
  
  (1) The Lyapunov exponents along all trajectories of the nonlinear
  flow $T_{c,n}^1$ are larger than $-\delta-\delta_0L> -\xi/k$ for
  every (fixed) $n\in U_\delta$ because the nonlinear perturbation has
  a Lipschitz constant less than $\delta_0L$ and all eigenvalues of
  the linear $\C^q$ flow $\dot E_c=H_c(n)E_c$ are larger than
  $-\delta$ for all $n\in U_\delta$.
%  More precisely, this means that,
%   uniformly for all $E_c\in\C^q$,
%   \begin{equation}\label{im:limsupt1cn}
%     \limsup_{t\to\infty}\left[\exp(-t\,\xi/k)
%     \left\|\left[\partial_2T_{c,n}^1(t;E_c)\right]^{-1}
%     \right\|\right]=0\mbox{.}
%   \end{equation}  
  
  (2) For any $n\in U_\delta$ there exists a simple open set
  $O(n)\subset \C^q$ that has a smooth boundary $\partial O(n)$ and
  contains the open ball $B_{R+1}(0)$ such that the flow $T_{c,n}^1$
  is overflowing invariant with respect to $O(n)$. That is,
  $T_{c,n}^1(1;O(n))\supset O(n)$ and the time $1$ image
  $T_{c,n}^1(1;\partial O(n))$ of the boundary $\partial O(n)$ has a
  positive distance from $\partial O(n)$.  This follows from the fact
  that the spectrum of $H_c(n)+\delta_0$ is positive for all $n\in
  U_\delta$ and $r_1=1$ for $\|E_c\|\geq R+1$.
  
  (3) We can find a smooth Lipschitz continuous function
  $r_2:U_\delta\times\C^q\mapsto[0,1]$ such that $r_2(n,E_c)=0$ if
  $E_c\in O(n)$ and $r_2(n,E_c)=1$ if $E_c\notin T_{c,n}^1(1;O(n))$.
  Due to the overflowing invariance (2) with respect to $O(n)$ this
  implies that the Lyapunov exponents of the flow
  $T_{c,n}^2:\R\times\C^q\mapsto\C^q$ generated by
  \begin{equation}\label{im:decdt}
    \dot E_c=(1-r_2(n,E_c))[H_c(n)+\delta_0r_1(\|E_c\|)]E_c
  \end{equation}
  are bigger than $-\xi/k$ along all trajectories of $T_{c,n}^2$.

  We introduce a smooth and globally Lipschitz continuous
  map $N:\R^m\to \R^m$ such that
  \begin{displaymath}
    N(n)=
    \begin{cases}
      n & n\in \mathcal{ N}_n\\
      \in U_\delta%\setminus\mathcal{N}_n 
      & \mbox{otherwise.}
    \end{cases}
  \end{displaymath}
  This is possible because $\mathcal{N}_n$ has a positive distance to
  the boundary of $U_\delta$. Furthermore, we can choose $N(n)$ such
  that $\partial_nN(n)=0$ if $\dist(n,\mathcal{N}_n)\geq1$. All
  statements of Theorem~\ref{im:thm} are concerned only with densities
  $n\in \mathcal{N}_n$, where $n=N(n)$. Thus, it is sufficient to
  prove Theorem~\ref{im:thm} for a modification of system
  \eqref{im:struc} that replaces $n$ by $N(n)$ in the right-hand-side
  of \eqref{im:struc}, which extends the definitions of the
  $n$-dependent quantities $P_c$, $P_s$, $B$, $\tilde P_c$, $H_c$,
  $r_2$ and $T_{c,n}^2$ $\R^m$ (by replacing the argument $n$ by
  $N(n)$).  We will abbreviate the notation by dropping the inserted
  modification of the argument $n$.
  
  The statements in \cite{BLZ99,BLZ00} about persistence and
  foliations assume the existence of a compact normally hyperbolic
  invariant manifold for a given unperturbed semiflow.  Our goal is to
  apply these theorems to the semiflow generated by system
  \eqref{im:struc} with $ \epsilon=0$.  For $\epsilon=0$ system
  \eqref{im:struc} possesses the invariant manifold
  \begin{displaymath}
    \mathcal{C}_0:=\{(E,n): n\in U_\xi, P_s(n)E=0\}\mbox{,}
  \end{displaymath}
  which is not compact. Thus, we will construct the modification of
  system~\eqref{im:struc} such that the modified system has a compact
  normally hyperbolic invariant manifold that is identical to
  $\mathcal{C}_0$ in the region of interest. We will achieve this by
  the insertion of $N(n)$ instead of $n$, and by a
  further modification of system~\eqref{im:struc} outside of
  $\mathcal{R}^{-1}\mathcal{N}$, combined with an extension by a
  one-dimensional auxiliary variable $x$. This extension will complete
  the cut-off version of $\mathcal{C}_0$ to an invariant sphere.
  
  Let $r_3:\R^m\mapsto[0,1]$ be a smooth function such that
  \begin{displaymath}
    r_3(n)=
    \begin{cases}
      0 & \mbox{if $\dist(n,\mathcal{N}_n)\leq1$,}\\
      1 & \mbox{if $\dist(n,\mathcal{N}_n)\geq2$.}
    \end{cases}
  \end{displaymath}
  Let $s:C^q\times\R^m\times\R\mapsto\R$ be a smooth function that is globally
  Lipschitz continuous along with its first derivative and satisfies
  {\renewcommand{\theenumi}{\alph{enumi}}%
    \begin{enumerate}
    \item\label{im:smx} $s(E_c,n,x)=(1-x^2)/2$ if
      $r_2(E_c,n)\neq 1$ and $|x|\leq2$,
    \item $\partial_2s(E_c,n,x)=0$ if $r_3(n)\neq1$,
    \item the set $\{(E_c,n,x): s(E_c,n,x)=0\}$ is a smooth compact manifold
      that is $C^\infty$ diffeomorphic to a $2q+m$ sphere in $\R^{2q+m+1}$ if we
      identify $\C^q$ with $\R^{2q}$, and
    \item\label{im:ds1} if $s(E_c,n,x)=0$ then the derivative of $s$ satisfies
      $\|\partial_1s(E_c,n,x)\|^2+ \|\partial_2s(E_c,n,x)\|^2+
      [\partial_3s(E_c,n,x)]^2\geq1$.
    \end{enumerate}
  }%
  Property~(\ref{im:smx})  implies that the manifold $\{s=0\}$ contains
  the set $\{(E_c,n,x): |x|=1,r_2(n,E_c)\neq1\}$.
  Property~(\ref{im:ds1}) clearly holds in this subset of $\{s=0\}$. A
  function $s$ satisfying (\ref{im:smx})--(\ref{im:ds1}) can be constructed by
  superimposing $s(E_c,n,x)=\|E_c\|^2+\|n\|^2+x^2-1$ with an appropriate
  ($\R$-)diffeomorphism of $C^q\times\R^m\times\R$ and then truncating $s$ for
  large arguments. 
  
  Consider the following modification (and extension) of
  system~\eqref{im:struc}:
  \begin{align}
    \allowdisplaybreaks
    \frac{d}{dt} E =\ & H(n)E -r_2(n,\tilde P_c(n)E)\left[H(n)P_c(n)E+
      \xi\, s\, B(n)\,
    \partial_1s^T \right] \label{im:te}\\
       & +\left[1-r_2(n,\tilde P_c(n)E)\right]\,\delta_0r_1(\|\tilde
    P_c(n)E\|)\,P_c(n)E\nonumber \\   
    \frac{d}{dt} n\ =\ &\epsilon F(N(n),E)- 
    r_3(n)\left[\epsilon F(N(n),E)^{\phantom{T}}+\xi s\,
    \partial_2s^T\right]\label{im:tn} \\
    \frac{d}{dt} x\  =\ &  -\xi\, s\,\partial_3s
    \label{im:tx}    
  \end{align}
  where $s$ and $\partial_js$ have been evaluated always at $(\tilde
  P_c(n)E,n,x)$. The right-hand-side of equation~\eqref{im:te} is
  differentiable with respect to $E$ if $\C$ is identified with
  $\R^2$.  In this sense $\partial_1s^T$ (the transpose of
  $\partial_1s$) contains $2q$ real components. System
  \eqref{im:te}--\eqref{im:tx} defines a semiflow $\tilde
  S(t;(E,n,x))$ on $X\times\R^m\times\R$ that is strongly continuous
  in $t$ and smooth with respect to $(E,n,x)$ because the
  right-hand-side is a smooth and globally Lipschitz continuous
  perturbation of a generator of a $C_0$ semigroup (note that the
  argument $n$ of $H$, $P_c$, $\tilde P_c$, and $r_2$ has been
  modified to $N(n)$; moreover all nonlinearities in the Nemytskii operators
  introduced by the modification have finite-dimensional ranges). Let
  us compare a trajectory $\tilde S(t;(E_0,n_0,x_0))=(\tilde
  E(t),\tilde n(t),\tilde x(t))$ generated by system
  \eqref{im:te}--\eqref{im:tx} to the trajectory
  $S(t;(E_0,n_0))=(E(t),n(t))$ generated by the original system
  \eqref{im:struc}. If $(E_0,n_0)\in\mathcal{R}^{-1}\mathcal{N}$ then
  $\tilde E(t)=E(t)$ and $\tilde n(t)=n(t)$ as long as $(E(t),n(t))$
  stays in $\mathcal{R}^{-1}\mathcal{N}$ because $r_1(\|\tilde
  P_c(n)E\|)=r_2(n,\tilde P_c(n)E)=r_3(n)=0$ if
  $(E,n)\in\mathcal{R}^{-1}\mathcal{N}$. In particular, the first two
  components of $\tilde S$, $\tilde E(t)$ and $\tilde n(t)$, do not
  depend on $x_0$ as long as $(\tilde E(t),\tilde n(t))$ stays in
  $\mathcal{R}^{-1}\mathcal{N}$.

  This relation between $\tilde S(t;\cdot)$ and $S(t;\cdot)$ enables us to
  prove the statements of Theorem~\ref{im:thm} by applying the results of
  \cite{BLZ98,BLZ99,BLZ00} to \eqref{im:te}--\eqref{im:tx} and ignoring the
  dynamics of the third component $x(t)$.

  \emph{Verification of conditions of \cite{BLZ98,BLZ99,BLZ00}}\\
  Consider the set
  \begin{displaymath}
    \tilde{\mathcal{C}}_0=\{(E,n,x)\in X\times\R^m\times\R:
    s(\tilde P_c(n)E,n,x)=0,P_s(n)E=0\}\mbox{.}
  \end{displaymath}
  The set $\tilde{\mathcal{C}}_0$ is the zero set of a smooth map. As
  it is non-empty (any $(E,n,x)$ with $E=0$, $n\in\mathcal{N}_n$,
  $x=1$ is in $\tilde{\mathcal{C}}_0$), it constitutes a $2q+m$
  dimensional manifold if the kernel of the linearization of the map
  in any element of $\tilde{\mathcal{C}}_0$ has dimension $2q+m$. The
  kernel of the linearization of the map $[s(\tilde P_c(n)E,n,x),
  P_s(n)E]^T$ in $(E_0,n_0,x_0)\in\tilde{\mathcal{C}}_0$ satisfies
  \begin{align}
    0&=\partial_1s\,E_c +\partial_2s\, n +\partial_1s\left[\partial_n\tilde
      P_c-\tilde P_c \partial_nP_s\right]nE_0+\partial_3s\,x\label{im:c0ds}\\
    0&=P_sE+P_s\partial_nP_snE_0\label{im:c0ps}
  \end{align}
  where we used $(E=B(n_0)E_c+P_s(n_0)E,n,x)$ as linearized variables,
  and the identities
  $\partial_nP_s=\partial_nP_s\,P_s+P_s\,\partial_nP_s$ and
  $P_sE_0=0$. Moreover, we always dropped the base point argument
  $(E_0,n_0,x_0)$ from $\partial_js$, $P_c$, $P_s$, and $\tilde P_c$.
  Equation~\eqref{im:c0ps} determines the $P_s(n_0)$ component of $E$.
  Thus, the kernel can have at most (real) dimension $2q+m+1$. Since
  $(\partial_1s,\partial_2s,\partial_3s)\neq 0$,
  equation~\eqref{im:c0ds} for $E_c$, $n$, and $x$ is non-trivial.
  Thus, the kernel has dimension $2q+m$. On the other hand, this is
  also the minimal dimension of the kernel since
  equation~\eqref{im:c0ps} has a zero $P_c$ component.  Consequently,
  $\tilde{\mathcal{C}}_0$ is a manifold of dimension $2q+m$. The
  manifold is bounded because the set
  $\{(E_c,n,x)\in\C^q\times\R^m\times\R:s(E_c,n,x)=0\}$ is bounded
  (diffeomorphic to a sphere), and $\|B(n)\|\|\tilde
  P_c(n)E\|+\|P_s(n)E\|$ provides an upper bound for $\|E\|$. Since
  $\tilde{\mathcal{C}}_0$ is finite-dimensional and closed (as a zero
  set of a continuous function) it is compact.

  The manifold $\tilde{\mathcal{C}}_0$ is invariant with
  respect to system \eqref{im:te}--\eqref{im:tx} at $\epsilon=0$ because
  \begin{align}
    \frac{d}{dt} s&=
    \partial_1s\cdot\left[\frac{d}{dt}[\tilde P_c(n)]E+
      \tilde P_c(n)\frac{d}{dt}E\right]+
    \partial_2s\cdot\frac{d}{dt}n+
    \partial_3s\cdot\frac{d}{dt}x\nonumber\\
    &=-\xi\, s\cdot\left[\partial_1s\cdot\partial_1s^T+
      \partial_2s\cdot\partial_2s^T+
      \partial_3s^2\right]\label{im:ds}\\
    \frac{d}{dt}\left[P_s(n)E\right]&=
    H(N(n))\,P_s(n)E\mbox{.}\label{im:des}
  \end{align}
  Here we dropped the argument $(\tilde P_c(n)E,n,x)$ of $s$ in
  \eqref{im:ds} and exploited that
  \begin{displaymath}
    0=\frac{d}{dt}[\tilde P_c(n(t))]=\frac{d}{dt}[P_c(n(t))]=
    \frac{d}{dt}[P_s(n(t))]    
  \end{displaymath}
  for \eqref{im:te}--\eqref{im:tx} at $\epsilon=0$ due to $\dot n=0$ (for
  $\dist(n,\mathcal{N}_n)\leq1$) or $\partial_nN=0$ (for
  $\dist(n,\mathcal{N}_n)\geq1$), keeping in mind the modification of
  the argument $n$ to $N(n)$. In fact, \eqref{im:ds}, \eqref{im:des}
  imply that the manifold $\tilde{\mathcal{C}}_0$ is
  exponentially attracting with rate $\xi$ (due to
  property (\ref{im:ds1}) of $s$ and  
  $\re\left[\spec H(N(n))\vert_{X_s(N(n))}\right]<-\xi$ for
  all $n\in\R^m$). 

  In order to check normal hyperbolicity of $\tilde{\mathcal{C}}_0$ we
  have to compute the linearization of $\tilde S(t;\cdot)$ along a
  trajectory $\tilde
  S(t;(E_0,n_0,x_0))=(E(t),n(t),x(t))\subset\tilde{\mathcal{C}}_0$.
  The trajectory satisfies
  \begin{displaymath}
    E(t)=
    \begin{cases}
      B(n_0)E_c(t) & \mbox{if $r_2(n_0,\tilde P_c(n_0)E_0)\neq1$,}\\
      E_0 &\mbox{otherwise,}
    \end{cases}\qquad\qquad
    n(t)=n_0\mbox{,}\quad x(t)=x_0
  \end{displaymath}
  where $E_c(t)$ satisfies (for $\epsilon=0$) the ODE~\eqref{im:decdt}
  with $n=n_0$ starting from $E_c(0)=\tilde P_c(n_0)E_0$. Thus,
  \begin{equation}\label{im:nhtan}
    \lim_{t\to\infty} \exp(t\xi/k)\left\|\partial_2\tilde
    S(t;(E_0,n_0,x_0))(E,n,x)\right\|=\infty
  \end{equation}
  uniformly for all $(E,n,x)$ in the tangent space
  $T\tilde{\mathcal{C}}_0$ of $\tilde{\mathcal{C}}_0$ with
  $\|(E,n,x)\|=1$. This follows from the corresponding property of
  \eqref{im:decdt} and the fact that $\partial_2\tilde S$, restricted
  to the tangent space of $\tilde{\mathcal{C}}_0$, is the identity if
  $r_2(n_0,\tilde P_c(n_0)E_0)=1$.

  For the transversal direction the equations \eqref{im:ds},
  \eqref{im:des} guarantee that
  \begin{equation}\label{im:nhnorm}
    \lim_{t\to\infty} \exp(t\xi)\left[
      |s(\tilde P_c(n(t))E(t),n(t),x(t))|+
      \|P_s(n(t))E(t)\|\right]=0
  \end{equation}
  for all trajectories $(E(t),n(t),x(t))$ of $\tilde S$. Since
  $\tilde{\mathcal{C}}_0$ is the regular zero set of $[s(\tilde
  P_c(n)E,n,x),P_s(n)E]^{\,T}$ this implies that the linearization of
  $\tilde{\mathcal{C}}_0$ has a transversal subbundle
  $N\tilde{\mathcal{C}}_0$ that is invariant under $\partial_2\tilde
  S$ with the same decay rate $\xi$. Hence, the growth conditions
  \eqref{im:nhtan} and \eqref{im:nhnorm} give the uniform normal
  hyperbolicity
  \begin{equation}\label{im:specgap}
    \left\|\partial_2\tilde S(t;(E_0,n_0,x_0))\vert_{N
        \tilde{\mathcal{C}}_0}\right\|<\exp(-t\xi)<
      \left\|\partial_2\tilde S(t;(E_0,n_0,x_0))\vert_{T
        \tilde{\mathcal{C}}_0}\right\|^{\,k}
  \end{equation}
  with spectral gap of size $k$
  for all sufficiently large $t>0$.

  \emph{Persistence of invariant manifold}\\
  The previous paragraph showed that the general theorems of
  \cite{BLZ98}, \cite{BLZ99}, \cite{BLZ00} can be applied to the
  invariant manifold $\tilde{\mathcal{C}}_0$ of the semiflow $\tilde
  S$, governed by system~\eqref{im:te}--\eqref{im:tx}, with
  $\epsilon=0$. According to \cite{BLZ99} this manifold persists under
  $C^1$-small perturbations such as a change to non-zero $\epsilon$.
  It stays $C^k$-smooth due to the spectral gap \eqref{im:specgap}.
  Thus, there exists an $\epsilon_0$ such that, for all
  $\epsilon\in[0,\epsilon_0)$, $\tilde S(t;\cdot)$ has a compact
  invariant $C^k$ manifold $\tilde{\mathcal{C}}$ which is a
  $C^1$-small perturbation of $\tilde{\mathcal{C}}_0$.  This implies
  that $P_s(n)E$ can be represented as a $C^k$-graph
  \begin{displaymath}
    P_s(n)E=\eta_0(\tilde P_c(n)E,n,x,\epsilon)    
  \end{displaymath}
  for $\tilde{\mathcal{C}}$. The evolution of $E$ and $n$
  does not depend on $x$ if $(\tilde P_c(n)E,n)\in\mathcal{N}$. Hence,
  $\eta_0(\tilde P_c(n)E,n,x,\epsilon)$ does not depend on $x$ if
  $(\tilde P_c(n)E,n)\in\mathcal{N}$. The range of $\tilde P_c(n)$ is
  $\C^q$, which allows us to parametrize the manifold
  $\tilde{\mathcal{C}}$ as a graph of $(E_c,n,\epsilon)\in
  \C^q\times\R^m\times[0,\epsilon_0)$. The semiflow $S$,
  governed by the original system~\eqref{im:struc}, and the $E$- and
  $n$-components of $\tilde S$, governed by the modification
  \eqref{im:te}--\eqref{im:tx}, are identical in $\mathcal{
    R}^{-1}\mathcal{N}$. Consequently, the set
  \begin{displaymath}
    \mathcal{C}:=\{( B(n)E_c+\eta_0(E_c,n,\epsilon),n):
    (E_c,n)\in \mathcal{N}\}
  \end{displaymath}
  is an invariant $C^k$ manifold of $S$ relative to $\mathcal{
    R}^{-1}\mathcal{N}$. By applying the projection $\tilde P_c(n)$ to
  the first equation of \eqref{im:struc} we obtain that the flow on
  $\mathcal{C}$ is governed by
  \begin{equation}
    \label{im:flowunexp}
    \begin{split}
      \frac{d}{dt} E_c &=\ \left[H_c(n)+
      \epsilon a_{1}(E_c,n,\epsilon)\right]E_c+ %\\ &\ +
      \epsilon a_{2}(E_c,n,\epsilon) \eta_0(E_c,n,\epsilon)\\ 
      \frac{d}{dt} n &=\ \epsilon F(n,B(n)E_c+\eta_0(E_c,n,\epsilon))
    \end{split}
  \end{equation}
  where the coefficients $a_1$ and $a_2$ are defined by
  \begin{displaymath}
    \begin{split}
      a_1(E_c,n,\epsilon)&=
      -\tilde P_c(n)\partial_nB(n)F(n,B(n)E_c+\eta_0(E_c,n,\epsilon))\\
      a_2(E_c,n,\epsilon) &=\ 
      B(n)^{-1}\partial_nP_c(n)F(n,B(n)E_c+\eta_0(E_c,n,\epsilon))
      P_s(n)\mbox{.}
    \end{split}
  \end{displaymath}

  \emph{Expansion of the graph $\eta_0$}\\
  The graph $\eta_0$ satisfies
  \begin{equation}
    \label{im:etaeps0}
    \eta_0(E_c,n,0)=0\mbox{\qquad for all
      $(E_c,n)\in \mathcal{N}$.}    
  \end{equation}
  Furthermore, the manifold $\mathcal{E}:=\{(E,n)\in X\times U_\delta:
  E=0\}$ is invariant with respect to $S$ for positive $\epsilon$. On
  $\mathcal{E}$, $\dot E=0$, and $\dot n=\epsilon F(n,0)$.
  Consequently, $\mathcal{E}\cap \mathcal{
      R}^{-1}\mathcal{N}\subset \mathcal{C}$, i.e.,
  \begin{equation}
    \label{im:etae0}
    \eta_0(0,n,\epsilon)=0\mbox{\qquad for $n\in \mathcal{N}_n$,
      $\epsilon\in[0,\epsilon_0)$.}
  \end{equation}
  Exploiting the fact that $\eta_0$ is $C^k$ with $k\geq2$, and the
  identities \eqref{im:etaeps0} and \eqref{im:etae0} we expand
  \begin{eqnarray}
    \eta_0(E_c,n,\epsilon) &=& \int_0^1\partial_1\eta_0(sE_c,n,\epsilon)\,ds E_c
    \nonumber \\
    &=&
    \epsilon\int_0^1\int_0^1\partial_1\partial_3\eta_0(sE_c,n,r\epsilon)\,dr\,ds
    E_c\mbox{.}
    \label{im:etaexp}
  \end{eqnarray}
  Denoting the double integral term in \eqref{im:etaexp} by $\nu$, we obtain
  \begin{equation}
    \label{im:nudef}
    \eta_0(E_c,n,\epsilon)=\epsilon \nu(E_c,n,\epsilon) E_c\mbox{.}
  \end{equation}
  We obtain the assertion \ref{im:flow} of the theorem by inserting
  \eqref{im:nudef} into system \eqref{im:flowunexp} for the flow on
  $\mathcal{C}$. The invariance of $\nu$ with respect to rotation of
  $E_c$ is a direct consequence of the rotational symmetry of the
  semiflow $S$.
  
  \emph{Exponential attraction toward $\mathcal{C}$}\\
  The theorems of \cite{BLZ98}, \cite{BLZ99}, \cite{BLZ00} imply that
  the set of all points that stay in a small tubular neighborhood of a
  compact normally hyperbolic invariant manifold $\mathcal{M}$ for all
  $t\geq 0$ form a center-stable manifold which is foliated by stable
  fibers of attraction rate close to the generalized Lyapunov numbers
  in the stable part of the linearization of the semiflow along
  $\mathcal{M}$. For $\tilde S$, governed by the modified
  system~\eqref{im:te}--\eqref{im:tx}, a whole tubular neighborhood
  $U$ of $\tilde{\mathcal{C}}$ is attracted by $\tilde{\mathcal{C}}$.
  (We can choose a uniform $U$ for all $\epsilon\in[0,\epsilon_0)$.)
  Thus, $U$ is foliated by stable fibers.
  
  For $\epsilon=0$, $P_s(n)E$ is exponentially decaying with rate $\xi$
  if $n\in\mathcal{N}_n$. If $\epsilon_0$ is sufficiently small
  there exists a $t_0\geq0$ such that $S(t_0;\Upsilon)\subset U$,
  which is foliated by stable fibers. 
  Hence, there exists a constant $M$ such that for all $u\in U$
  there exists a fiber base point $u^*\in\tilde{\mathcal{ C}}$ such
  that
  \begin{equation}\label{im:attrineq}
    \|\tilde S(t;u)-
    \tilde S(t;u_*)\| 
    \leq Me^{-\Delta t}
  \end{equation}
  where we may have to decrease $\epsilon_0$ (if necessary) in order to keep
  the decay rate at $\Delta$ in \eqref{im:attrineq}.
  
  Let $t_1\geq 0$ be such that $Me^{-\Delta t_1}$ is less than the
  distance between the set $\mathcal{ R}\Upsilon$ and the boundary of
  $\mathcal{ N}$. Then, we can choose $t_c=t_0+t_1$ to obtain
  assertion~\ref{im:attr} of the theorem: Let $(E,n)\in\Upsilon$ and
  $t\geq0$ be such that $S([0,t+t_c];(E,n))\subset\Upsilon$. Then the
  $E$- and the $n$-component of $u=\tilde S(t_c;(E,n,1))$ are in $U$,
  and, furthermore, the fiber base point $u_*=(E_*,n_*,x_*)$ for $u$
  satisfies $(E_*,n_*)\in\mathcal{C}$ and $x_*=1$. Since $\tilde S$
  (ignoring the evolution of the auxiliary variable $x$) and $S$ are
  identical, inequality \eqref{im:attrineq} implies the inequality
  \eqref{im:attreq} for $(E_*,n_*)$.\qquad \qed\end{pf}

\section{Practical applications and conclusions}
\label{sec:out}

\begin{figure}[ht!]
  \centering
  \includegraphics[scale=0.5]{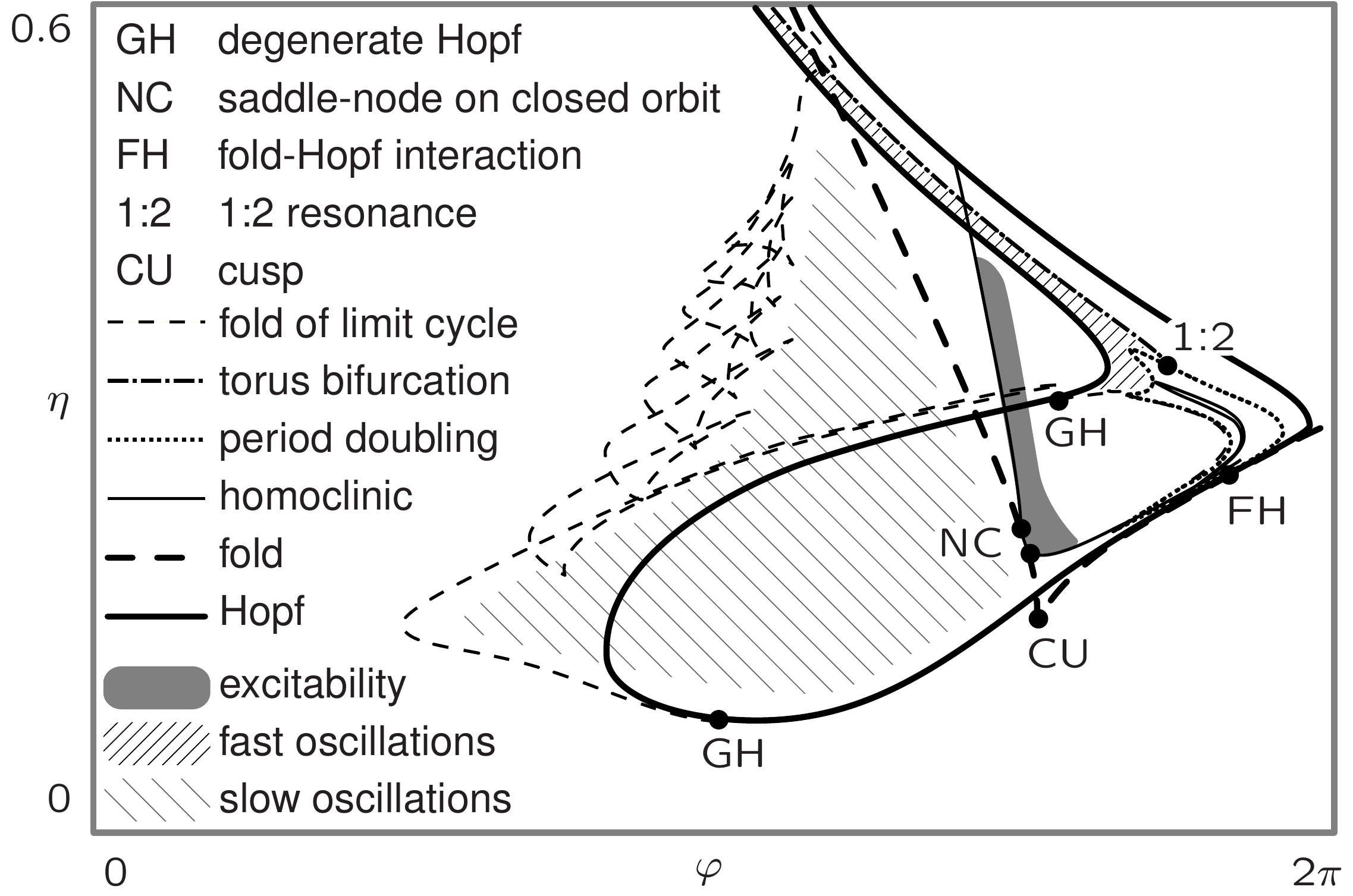}
  \caption{Bifurcation diagram for the two-section laser investigated in \cite{WBRH02}. 
    The parameters are: $l_2=1.136$, $r_0=10^{-5}$, $r_L=\eta
    e^{i\phi}$, $d_1=-0.275$, $\kappa_1=3.96$, $\tilde g_1=2.145$
    (linear gain model), $\alpha_1=5$, $\rho_1=0.44$, $\Gamma_1=90$,
    $\Omega_{r,1}=-20$, $I_1=6.757\cdot10^{-3}$, $\tau_1=359$,
    $\dot n_2=\kappa_2=\beta_2=\rho_2=0$. The bifurcation parameters are the
    strength $\eta$ and the phase $\phi$ of the feedback from the
    facet $r_L$ of section $S_2$. In the experiment these parameters
    can be varied by changing the current in $S_2$. The highlighted
    dynamical regimes are of particular practical interest.}
  \label{fig:bif}
\end{figure}

\paragraph*{Truncated reduction}
The graph of the invariant manifold enters the description
\eqref{im:flowc} of the flow on $\mathcal{C}$ only in the form
$O(\epsilon^2)\nu$.  The consideration of system \eqref{im:flowc}
(truncating $\nu$ to $0$) has been extremely useful for numerical and
analytical investigations of dynamics of multi-section semiconductor
lasers because the dimension of system \eqref{im:flowc} is typically
low ($q$ is often either $1$ or $2$); see, e.g.,
\cite{WBRH02,B94,BRS98,BWSM97,S01a,WBWR95,UBBWH04}. For illustration,
Fig.~\ref{fig:bif} shows a two-parameter bifurcation diagram for a
two-section laser \cite{WBRH02}. After reduction of the rotational
symmetry the dimension of the invariant manifold $\mathcal{ C}$ is $4$
($m=1$ since section $\dot n_2=0$, and $q=2$) in the parameter range
covered by the diagram. A detailed numerical comparison of
Fig.~\ref{fig:bif} with simulation results for the PDE model
\eqref{twe:twe}--\eqref{twe:carrier} and more accurate models can be
found in \cite{RW02}.

\paragraph*{Delay-differential equations}
Theorem \ref{im:thm} also applies to systems of delay-differential
equations (DDEs) as they are widely used in laser dynamics (such as,
for example, the Lang-Kobayashi system for delayed optical feedback or
two mutually coupled lasers; see \cite{TA98,KL00} and references
therein).  They also have the structure of system~\eqref{int:geneq}
where $E\in C([-1,0];\C)$, and $H(n)$ is a linear delay operator.
% \begin{equation}
%   \label{out:lke}
%   \begin{split}
%     \frac{d}{dt}E(t) &=\ (1+i\alpha)nE(t)+\eta e^{i\phi} E(t-1)\\
%     \frac{d}{dt}n(t) &=\ \epsilon\left(F(n)-g(n)|E(t)|^2\right)
%   \end{split}
% \end{equation}
% if its scaling is appropriate to the situation of a short external
% cavity \cite{WT02p}.  System \eqref{out:lke} generates a semiflow in
% the Banach space $C([-1,0];\C)\times\R$ and has also the structure
% \eqref{int:geneq}. The parameters have the same sense as in
% \eqref{twe:twe}--\eqref{twe:carrier} (we have dropped the indices
% since there is only one section).
The parameter $\epsilon$ is small if the feedback cavity is short.
According to \cite{DGLW95}, $H(n)$ generates an eventually compact
semigroup and, thus, the existence of critical carrier densities with
a spectral gap can be shown analytically \cite{KL00}. Moreover, the
cut-off modification performed in the proof of Theorem \ref{im:thm}
manipulates only the finite-dimensional components $E_c$ and $n$.
Hence, the proof does not rely on the ability to cut-off a smooth map
smoothly in the infinite-dimensional space $X$ which is the Hilbert
space $X=\Lint^2([0,L];\C^2)\times\Lint^2([0,L];\C^2)$ in
Section~\ref{sec:im} but a Banach space for DDEs.  The only property
of the operator $H(n)$ used in the proof is the existence of a
spectral splitting and the smooth dependence of the dominating
subspace $X_c$ on $n$.  Consequently, Theorem~\ref{im:thm} applies to
\eqref{int:geneq} if $H(n)$ is a delay operator, reducing the DDEs to
low-dimensional systems of ODEs.

% There are other models in the spirit of \eqref{out:lke} for different
% experimental situations, e.g., for lasers subject to dispersive
% feedback or for two lasers interacting with each other. All have the
% structure of \eqref{int:geneq} where $H$ is a delay operator smoothly
% depending on $n$, and $\epsilon$ is small if the external cavity is
% short. Hence, Theorem \ref{im:thm} allows to reduce these models
% locally to low-dimensional systems of ODEs.

\section*{Acknowledgments}
The research of J.S.\ was partially supported by the the Collaborative
Research Center 555 ``Complex Nonlinear Processes'' of the Deutsche
Forschungsgemeinschaft (DFG), and by EPSRC grant GR/R72020/01. The
author thanks Mark Lichtner and Bernd Krauskopf for discussions and
their helpful suggestions.

\bibliography{lasbifref}
\bibliographystyle{elsart-num}

\end{document}